\documentclass[10pt,twoside]{article}
\usepackage[latin1]{inputenc}
\usepackage{amsmath}
\usepackage{graphicx}
\usepackage{yhmath}
\usepackage[table]{xcolor}
\usepackage{ mathrsfs} % permet d'avoir le D cursive
\usepackage{amssymb}
\usepackage{url}

\input epsf
\def\limiten{\renewcommand{\arraystretch}{0.5}
\begin{array}[t]{c}\stackrel{}{\longrightarrow} \\
{\scriptstyle n\rightarrow
\infty}\end{array}\renewcommand{\arraystretch}{1}}

\def\limitepsn{\renewcommand{\arraystretch}{0.5}
\begin{array}[t]{c}\stackrel{a.s.}{\longrightarrow} \\
{\scriptstyle n \rightarrow
\infty}\end{array}\renewcommand{\arraystretch}{1}}

\def\limiteloin{\renewcommand{\arraystretch}{0.5}
\begin{array}[t]{c}\stackrel{{\cal D}}{\longrightarrow} \\
{\scriptstyle n\rightarrow
\infty}\end{array}\renewcommand{\arraystretch}{1}}

\def\limiteproban{\renewcommand{\arraystretch}{0.5}
\begin{array}[t]{c}\stackrel{{\cal P}}{\longrightarrow} \\
{\scriptstyle n\rightarrow
\infty}\end{array}\renewcommand{\arraystretch}{1}}

\numberwithin{equation}{section}

\newtheorem{thm}{Theorem}[section]

\newtheorem{lem}[thm]{Lemma}

\newcommand{\E}{\ensuremath{\mathbb{E}}}
\newcommand{\R}{\ensuremath{\mathbb{R}}}
\newcommand{\Z}{\ensuremath{\mathbb{Z}}}

\newcommand{\N}{\ensuremath{\mathbb{N}}}

\definecolor{grisclair}{gray}{0.9}

\renewcommand{\arraystretch}{.8}

\begin{document}
%\date{ }
\title{\bf Testing for parameter change in general integer-valued time series}
 \maketitle \vspace{-1.0cm}
 \begin{center}
   Mamadou Lamine DIOP$^{\text{a}}$ \footnote{Supported by AIRES-SUD (Appuis Intégrés pour le Renforcement des Equipes Scientifiques du Sud)} and
   William KENGNE$^{\text{b}}$ \footnote{Developed within the MME-DII center of excellence (ANR-11-LABEX-0023-01)  \url{http://labex-mme-dii.u-cergy.fr/} }.
 \end{center}

  \begin{center}
 {\it $^{\text{a}}$  LERSTAD, Université Gaston Berger, Saint-Louis, Sénégal.\\
 $^{\text{b}}$ THEMA, Université de Cergy-Pontoise, 33 Boulevard du Port, 95011 Cergy-Pontoise Cedex, France.\\
  E-mail: diopml@yahoo.fr; william.kengne@u-cergy.fr}
 \end{center}

 \pagestyle{myheadings}
 \markboth{Testing for parameter change in general integer-valued time series}{M.L. Diop and W. Kengne}

\textbf{Abstract} :
We consider the structural change in a class of discrete valued time series that the conditional distribution follows a one-parameter exponential family.
 We propose a change-point test based on the maximum likelihood estimator of the parameter of the model. Under the null hypothesis (of no change), the test statistics converges
 to a well known distribution, allowing for the calculation of the critical values of the test.
 The test statistic diverges to infinity under the alternative, that is, the test asymptotically has power one.
 Some simulation results and real data applications are reported to show the applicability of the test procedure.\\


\section{Introduction}

  In recent years, modeling time series of counts has become an important topic of statistical research, mainly because of the wide area of applications
  (epidemiology, economic, finance, actuarial science, etc), see for instance the of books of Cameron and Trivedi (2013) and Kedem  and Fokianos (2005) and the references therein.
   Several models have been proposed to describe time series count data. Those models lead to capture some specific phenomena
   (such as excess zeros, overdispersion, heteroscedasticity, etc) that often display counts data and also lead to describe the dependence structure of the observations. \\
  Let $Y=(Y_t)_{t \in \Z}$ be an integer-valued  time series ; let
   $\mathcal{F}_t$ be the filtration that represents all the information that is known at time $t$ (see below).
   In the setting of autoregressive models, there is a huge number of papers that focussed on the distribution of $Y_t/\mathcal{F}_{t-1}$, denoted by $p(\cdot/\mathcal{F}_{t-1})$. \\
    \indent  Neumann (2011) has assumed that
   $ Y_t/\mathcal{F}_{t-1}   \sim\mbox{ Poisson}(\lambda_t)$ with $\lambda_t = f(\lambda_{t-1},Y_{t-1})$ (where $f$ is a non-negative mesurable function) and focused
   on the stationarity and ergodicity of the bivariate process $(Y_t,\lambda_t)_{t \in \N}$. The models with linear and log-linear link function have been studied
   by Fokianos \textit{et al.} (2009) and  Fokianos and Tjøstheim (2011).
   We refer to Ferland \textit{et al.} (2006), Weiß (2009), Fokianos and Tjøstheim (2012), Doukhan \textit{et al.} (2012),  Zhu (2012a), Fokianos  and  Neumann (2013),
    Fried \textit{et al.} (2015),  etc  for some recent works on Poisson autoregressive models.
   Moysiadis and Fokianos (2014) studied the problem of ergodicity, stationarity in binary and categorical time series models with feedback, whereas
   Christou And Fokianos (2014) focussed on inference and diagnostics in negative binomial processes for count time series.
    Also, see Davis and Wu (2009) and Zhu (2011, 2012b) for model with negative binomial and a distribution with zero-inflated characteristic for modeling data with
     excess zeros.\\
   Recently, Davis and Liu (2012)  focussed on the model where the conditional distribution $Y_t/\mathcal{F}_{t-1}$ belongs to the class of
   one-parameter exponential family ; which contains Poisson and negative binomial (with fixed number of failures) distribution.
    They established the stationarity and the absolute regularity properties of the process as well as the
     consistency and  asymptotic normality of the maximum likelihood estimator of the parameter of the model.

   On the other hand, as it has been pointed out in numerous works, many real data often exhibit structural change occurred during the data collecting processes.
   Ignoring these breaks can adversely affect any statistical inference on such data; as suggested by Berkes {\it et al.} (2006), Kengne (2012), Franke \textit{et al.} (2012),
   Doukhan and Kengne (2015), just to name a few.
   Change-point detection is now an important field in time series analysis ; in view of the huge number of papers written in this direction during the last three decades.
   Two approaches are usually considered : the on-line and off-line detection.
   On-line detection focussed on sequential change detection as and when new data arrive whereas off-line approach leading to change-point detection when all data are available.
   See Basseville  and Nikiforov (1993) and Cs{\"o}rg{\"o} and Horv{\'a}th (1997) for surveys on these approaches.

    We focuss here on the off-line change-point detection for time series of counts models. This problem has attracted some attention in the recent past.
   Kang and Lee (2009) proposed the cumulative sum (CUSUM)  procedure for detecting changes in random coefficient
   integer-valued autoregressive models with Poisson innovations. Franke \textit{et al.} \cite{Franke_2012} focussed on the residuals CUSUM procedures for testing for parameter change
   in Poisson autoregressive models. See Hudecov{\'a} (2013) and  Fokianos {\it et al.} (2013) for change-point detection in binary time series.
    Kang and Lee (2014) proposed procedures based on the conditional maximum likelihood estimator for change-point detection in the one order Poisson autoregressive models.
  % See also Fokianos and Fried (2010, 2012) for mean shift detection in the linear and log-linear Poisson autoregression.
   In the same vein, Doukhan and Kengne \cite{Doukhan_CH_P} proposed  two tests  based on the likelihood of the observations for change-point detection in a general class
   of Poisson autoregressive models.

In this contribution, we consider a process $\{Y_{t},t\in \mathbb{Z}\}$ satisfying :
\begin{equation} \label{Model}
Y_t|\mathcal{F}_{t-1}\sim p(y|\eta_t)~~\textrm{with}~~X_t:=\mathbb{E}(Y_t|\mathcal{F}_{t-1})=f_{\theta_0}(X_{t-1},Y_{t-1})~~\textrm{and}~~X_{t}=A^{\prime}(\eta_t)
\end{equation}
where $\mathcal{F}_{t-1}=\sigma\left\{\eta_1, X_{t-1},X_{t-2},\cdots \right\}$ is the $\sigma$-field generated by the whole past at time $t-1$,
  $p(.|.)$ is a discrete distribution that followed a one-parameter exponential family ; that is
\[ p(y|\eta)=\exp\left\{\eta y -A(\eta)\right\}h(y) \]
where $\eta$ is the natural parameter, $A(\eta)$, $h(y)$ are known fonctions and $f_{\theta_0}(.)$ is a non-negative bivariate function defined on
$[0,+\infty)\times \mathbb{N}_0$ (where $\mathbb{N}_0 = \{0,1,2,\cdots \}$), assumed to be know up to a parameter $\theta_0 \in \Theta$; where $\Theta$ is a compact set of
 $ \mathbb{R}^{d}$ ($d \in \N$).
This distribution family contains many of the most common discrete distributions,
including Poisson, negative binomial (with fixed number of failures), Bernoulli, etc. \\

 This class of model have been studied by Davis and Liu \cite{Davis_2012}. They established the stationarity and the absolute regularity properties of the process.
The consistency and asymptotic normality of the maximum likelihood estimator of the parameter of the model are also proved.\\

Assume that a trajectory $(Y_1,\cdots,Y_n)$ of $\{Y_{t},t\in \mathbb{Z}\}$ are observed and consider the following hypothesis :
\begin{enumerate}
    \item [H$_0$:] The observations $Y_1,\cdots,Y_n$  are a trajectory of a process $\left(Y_t\right)_{t\in \mathbb{Z}}$ solution of (\ref{Model}), depending on $\theta_0 \in \Theta$.
    \item [H$_1$:] There exists $((\theta^{*}_1,\theta^{*}_2),t^{*}) \in \Theta^{2}\times \{2,3,\cdots, n-1 \}$ such that
     $(Y_1,\cdots,Y_{t^{*}})$ is a trajectory of a process $\{Y^{(1)}_t, t \in \mathbb{Z}\}$   and  $(Y_{t^{*}+1},\cdots,Y_n)$ a trajectory of
      $\{Y^{(2)}_t, t \in \mathbb{Z}\}$, where the processes $\{Y^{(1)}_t, t \in \mathbb{Z}\}$ and $\{Y^{(2)}_t, t \in \mathbb{Z}\}$ are stationary solutions of (\ref{Model})
      depending respectively on $\theta^{*}_1$ and $\theta^{*}_2$ with $\theta^{*}_1 \neq \theta^{*}_2$.
\end{enumerate}
 The change-point problem considered here is more general than those studied before, since the class of distribution considered in (\ref{Model}) contains among others,
 the Poisson, Bernoulli and negative binomial (with fixed number of failures) distribution.
 We will generalize the procedure proposed by Doukhan and Kengne \cite{Doukhan_CH_P} to the class of model (\ref{Model}) ; a usefulness of this generalization is the application
 to binary time series and negative binomial INGARCH models (see the simulation study and real data application).

  \medskip

 In Section 2, the assumptions and the definition of the likelihood estimator as well as some examples of the model (\ref{Model}) are provided.
 Section 3 is devoted to the procedure for change-point detection in the model (\ref{Model}). In Section 4, we conduct a simulation study.
 Applications to the number of transactions for the stock Ericsson B and the US recession data are presented in Section 5. The proofs of the main results
are provided in Section 6.

\section{Assumptions, likelihood inference and examples}\label{presentation}

 \subsection{Assumptions}

 Throughout this paper we will assume the classical Lipschitz-type condition on the model (\ref{Model}) :

 \medskip

\textbf{Assumption} (\textbf{A}($\Theta$)):
There exists two non-negative real numbers $\delta_{1}$ and $\delta_{2}$ satisfying $\delta_{1} +\delta_{2} <1$ and such that for any
 $(x, y), (x', y') \in [0,+\infty)\times \mathbb{N}_{0}$,
 \[ \sup_{\theta \in \Theta  }\left|f_{\theta}(x,y)-f_{\theta}(x',y')\right|\leq \delta_{1}\left|x-x'\right|+\delta_{2}\left|y-y'\right|. \]

Under the assumption (\textbf{A}($\Theta$)), Davis and Liu \cite{Davis_2012} have shown that the process $\{Y_{t},t\in \mathbb{Z}\}$  is absolutely regular
with geometrically decaying coefficients and $\left\{(X_{t},Y_{t}),~t\geq 1\right\}$ is strictly stationary and ergodic.
 In addition, the conditional mean can be expressed as a function of only the past observations, meaning,
 there exists  $f^{\theta}_\infty:[0,\infty)^{\infty}\rightarrow [0,\infty) $  such that
 \begin{equation}\label{f_inf}
  X_t=f^{\theta}_\infty(Y_{t-1},Y_{t-2},\cdots) ~ a.s. ~ .
 \end{equation}
These results generalize those obtained by Fokianos \textit{et al.} (2009), Fokianos and Tjøstheim (2012) and Neumann (2011) in the sense that the class of the conditional
distribution considered here is more large.

 \subsection{Likelihood inference}\label{Likelihood}
In this section, we make a brief overview of the maximum likelihood estimation in the model (\ref{Model}) with the main asymptotic properties.
Let $(Y_{1},\ldots,Y_{n})$ be a trajectory generated from the model (\ref{Model}), according to the true parameter $\theta_0$.
The likelihood function conditioned on $\eta_1$ is given for any $\theta \in \Theta$ by
$$\mathcal{L}(\theta | Y_{1},\ldots,Y_{n},\eta_1):=\prod_{t=1}^{n}\exp\left\{\eta_t(\theta )Y_{t}-A(\eta_t(\theta ))\right\}h(Y_{t})$$
where $\eta_t(\theta)=(A^{\prime})^{-1}(X_t(\theta))$ is updated through the relation
 \[X_t(\theta) := f_{\theta}(X_{t-1},Y_{t-1}). \]
The conditional log-likelihood function, up to a constant independent of $\theta$, is given by
\[ L_n(\theta) := \log\big(\mathcal{L}(\theta | Y_{1},\ldots,Y_{n},\eta_1)\big) = \sum_{t=1}^{n}\ell_t(\theta) ~~\textrm{with}~~\ell_t(\theta)= \eta_t(\theta )Y_{t}-A(\eta_t(\theta )). \]
The maximum likelihood estimator (MLE) of $\theta_0$ is defined by
 \[ \widehat{\theta}_{n}: =  \underset{\theta \in \Theta}{\text{argmax}} (L_n(\theta)).\]
To ensure the consistency and the asymptotic normality of the MLE, we impose the following regularity conditions (see also Davis and Liu \cite{Davis_2012}) :
\begin{enumerate}
    \item [(\textbf{A0}):] $\theta_0$ is an interior point in the compact parameter space $\Theta \subset \mathbb{R}^{d}$.
    \item [(\textbf{A1}):] For any $\theta \in \Theta$, $f^{\theta}_{\infty} \geq x^{*}_{\theta} \in \mathcal{R}(A^{\prime})$, where $\mathcal{R}(A^{\prime})$ is the range of $A^{\prime}(\eta)$. Moreover $x^{*}_{\theta}\geq x^{*} \in \mathcal{R}(A^{\prime})$ for all $\theta$.
    \item [(\textbf{A2}):] For any $y \in \mathbb{N}^{\infty}_{0}$,
                                   % $y \in [0,\infty)^{\infty}\cup \mathbb{N}^{\infty}_{0}$,
                         the mapping $\theta \mapsto f^{\theta}_{\infty}(y)$ is continuous.
    \item [(\textbf{A3}):] $\mathbb{E}\big[ Y_1 \big| \sup_{\theta \in \Theta} (A^{\prime})^{-1}(f^{\theta}_{\infty} (Y_0,Y_1,\cdots)) \big| \big] <\infty$.
    \item [(\textbf{A4}):] If there exists a $t \geq 1$ such that $X_t(\theta) = X_t(\theta_0)$ a.s., then $\theta=\theta_0$.
    \item [(\textbf{A5}):] The mapping $\theta \mapsto f^{\theta}_{\infty}$ is twice continuously differentiable.
     \item [(\textbf{A6}):] $\E \Big[ \sup_{\theta \in \Theta} \Big| A^{\prime \prime}(\eta_t(\theta)) \Big(\frac{\partial}{ \partial \theta_i}\eta_t(\theta)\Big) \Big(\frac{\partial}{ \partial \theta_j}\eta_t(\theta) \Big)  \Big| \Big] <\infty$
           ~and~ $ \E \Big[ \sup_{\theta \in \Theta} \Big| \big(Y_{t}-A^{\prime}(\eta_t(\theta ))\big) \frac{\partial^{2}} {\partial\theta_i \partial\theta_j} \eta_t(\theta)  \Big| \Big] < + \infty $
           for any $i,j = 1,\cdots, d$.
\end{enumerate}
Under the above assumptions and the Lipschitz-type condition (\textbf{A}($\Theta$)), Davis and Liu \cite{Davis_2012} have used the Lemma 3.11 in Pfanzagl (1969)
to show that the MLE $\widehat{\theta}_{n}$ is strongly consistent (under H$_0$), that is, $$\widehat{\theta}_{n}\stackrel{a.s.}{\longrightarrow}\theta_{0} ~~ \textrm{as}~~n\rightarrow \infty.$$
They have also proved the asymptotic normality (under H$_0$) of the MLE $\widehat{\theta}_{n}$, that is,
 \[\sqrt{n}(\widehat{\theta}_{n}-\theta_{0}) \stackrel{\mathcal{D}}{\longrightarrow}\mathcal{N}(0,\Omega^{-1}) ~~ \textrm{as}~~n\rightarrow \infty \]
where
\begin{equation}\label{def_Sigma}
\Omega=\mathbb{E} \big[ A^{\prime \prime}\left(\eta_{0}(\theta_0)\right)\left(\frac{\partial \eta_{0}(\theta_0)}{\partial \theta}\right)\left(\frac{\partial \eta_{0}(\theta_0)}{\partial \theta}\right)^{T} \big].
\end{equation}
Under H$_0$, a consistent estimator of $\Omega$ is (see Lemma \ref{Lem_H0})
 \[  \widehat{\Omega}_{n}=\left.\frac{1}{n}\sum_{t=1}^{n}\left({A}^{\prime \prime}\left( \eta_t(\theta)\right)\left(\frac{\partial \eta_t(\theta)}{\partial \theta}\right)\left(\frac{\partial \eta_t(\theta)}{\partial \theta}\right)^{T}\right)\right|_{\theta=\widehat{\theta}_n } .\]
 \subsection{Examples}

\subsubsection{Linear autoregressive models}
\begin{enumerate}
    \item   Consider an  INGARCH$(1,1)$ process $\left(Y_t\right)_{t\in \mathbb{Z}}$ defined by :
    \begin{eqnarray}\label{INGARCH(1,1)}
    Y_{t}|\mathcal{F}_{t-1}\sim\textrm{Poisson}(X_{t}), \text{ with } X_{t} = \alpha_0^* +\alpha^* Y_{t-1} + \beta^* X_{t-1}  ;
    \end{eqnarray}
     the true parameter $\theta_0 = (\alpha_0^*, \alpha^* , \beta^*) \in \Theta$ where $\Theta$ is a compact subset of $(0, +\infty) \times [0, +\infty)^2$
     such that $\alpha + \beta < 1$ for all $\theta = (\alpha_0, \alpha, \beta) \in \Theta$ (see  Ferland \textit{et al.} (2006)).
     This is a particular case of the model (\ref{Model}) with  $f_\theta(x,y)=\alpha_{0} +\alpha y + \beta x$, $\eta_t(\theta)=\log(X_t(\theta))$ and $A(\eta_t(\theta))=e^{\eta_t(\theta)}$
     for all  $\theta=\left(\alpha_{0}, \alpha, \beta \right) \in \Theta$.
    The assumption (\textbf{A}($\Theta$)) holds. This ensures the existence of a strictly stationary solution for the model (\ref{INGARCH(1,1)}) and the ergodicity of the process $\left\{(X_{t},Y_t)\right\}$.
    Moreover, by the recursive substitution, we get for all $\theta \in \Theta$,
    \[ X_{t}(\theta)=\frac{\alpha_{0}}{1-\beta} + \alpha \sum_{k=0}^{\infty}\beta^{k}Y_{t-1-k}:=f^{\theta}_\infty(Y_{t-1},\cdots). \]
 Hence, it is easily to see that (\textbf{A2}), (\textbf{A4}) and (\textbf{A5}) are satisfied.
 To check the conditions (\textbf{A1}), (\textbf{A3}) and (\textbf{A6}), suppose that the true parameter vector $\theta_0$ lies in a compact set :
\begin{eqnarray}\label{Set_Theta}
\Theta=\Big\{\theta=(\alpha_{0},\alpha,\beta) \in \mathbb{R}^{3}_{+}:~0<\alpha_{L}\leq \alpha_{0} \leq \alpha_{U},~\epsilon \leq \alpha + \beta \leq 1- \epsilon \Big\}~\rm{for~ some}~ \alpha_L, \alpha_U, \epsilon>0.
\end{eqnarray}
Since $Y_{t}\geq 0$ for all $t \in \mathbb{Z}$, we have $f^{\theta}_\infty(Y_{t-1},\cdots) \geq \frac{\alpha_{0}}{1-\beta} \geq \dfrac{\alpha_U}{\epsilon}$.
 Hence, (\textbf{A1}) holds.
 For the condition (\textbf{A3}), remark that for all $\theta \in \Theta$,
 $  \dfrac{\alpha_U}{\epsilon}  \leq f^{\theta}_\infty(Y_{t-1},\cdots) \leq  \dfrac{\alpha_{U}}{\epsilon} +\sum_{k=1}^{\infty}(1-\epsilon)^{k}Y_{1-k}   $.
 Therefore,
\begin{eqnarray*}
 \mathbb{E}\big[ Y_1 \big| \sup_{\theta \in \Theta} (A^{\prime})^{-1}(f^{\theta}_{\infty} (Y_0,Y_1,\cdots)) \big| \big] &\leq&
 \E\Big\{Y_1 \Big[ \big|\log(\frac{\alpha_{U}}{\epsilon})| +  | \log\Big(\frac{\alpha_{U}}{\epsilon} +\sum_{k=1}^{\infty}(1-\epsilon)^{k}Y_{1-k}\Big) \big| \Big] \Big\}\\
 &\leq& \E\Big\{Y_1 \Big[ 2 \big|\log(\frac{\alpha_{U}}{\epsilon}) \big| +  \frac{\alpha_{U}}{\epsilon} +\sum_{k=1}^{\infty}(1-\epsilon)^{k}Y_{1-k} \Big] \Big\}\\
&\leq&( 2 \big|\log(\frac{\alpha_{U}}{\epsilon}) \big| +  \frac{\alpha_{U}}{\epsilon} )\mu  +\sum_{k=1}^{\infty}(1-\epsilon)^{k}\left(\gamma_{Y}(k)+\mu^{2}\right)<\infty
\end{eqnarray*}
where $\mu=\E(Y_t)=\frac{\alpha^*_0}{1-\alpha^*-\beta^*}$ and $\gamma_{Y}(k)=\mu C(\theta_0)(\alpha^*+\beta^*)^{k-1}$ with $C(\theta_0)$ is a positive constant depending on $\theta_0$ (see \cite{Fokianos_2009}).
Hence, the condition (\textbf{A3}) is obtained.\\
 Furthermore, for $i=1,2,3$, we have $\frac{\partial}{\partial \theta_i} \eta_t(\theta)=\frac{1}{X_t(\theta)} \frac{\partial}{\partial \theta_i} X_t(\theta)$. Hence
$$A^{\prime \prime}(\eta_t(\theta)) \left(\frac{\partial}{ \partial \theta_i}\eta_t(\theta)\right) \left(\frac{\partial}{ \partial \theta_j}\eta_t(\theta)\right)  = \frac{1}{X_t(\theta)} \left(\frac{\partial}{\partial \theta_i} X_t(\theta)\right) \left(\frac{\partial}{\partial \theta_j} X_t(\theta)\right).$$
Therefore
\begin{eqnarray*}
\left|A^{\prime \prime}(\eta_t(\theta)) \left(\frac{\partial}{ \partial \theta_i}\eta_t(\theta)\right) \left(\frac{\partial}{ \partial \theta_j}\eta_t(\theta)\right) \right| \leq C \sup_{\theta \in \Theta}\left|\frac{\partial}{\partial \theta_i} X_t(\theta) \right| \cdot \sup_{\theta \in \Theta}\left| \frac{\partial}{\partial \theta_j} X_t(\theta) \right| = C \left\| \frac{\partial}{\partial \theta_i} f^{\theta}_\infty(Y_{t-1} \right\| _\Theta \left\| \frac{\partial}{\partial \theta_j} f^{\theta}_\infty(Y_{t-1}\right\| _\Theta,
\end{eqnarray*}
 with  $C= 1/\alpha_L$.
 Thus,
\begin{eqnarray}\label{Major1}
\E \left[ \sup_{\theta \in \Theta} \left|A^{\prime \prime}(\eta_t(\theta)) \left(\frac{\partial}{ \partial \theta_i}\eta_t(\theta)\right) \left(\frac{\partial}{ \partial \theta_j}\eta_t(\theta)\right) \right| \right] &\leq& C \cdot \E \left[\left\| \frac{\partial}{\partial \theta_i} X_t(\theta)\right\| _\Theta \left\| \frac{\partial}{\partial \theta_j} X_t(\theta)\right\| _\Theta \right] \nonumber\\
 &\leq& C \left\|\left\| \frac{\partial}{\partial \theta_i} X_t(\theta)\right\| _\Theta\right\|_2 \cdot \left\|\left\| \frac{\partial}{\partial \theta_j} X_t(\theta)\right\| _\Theta\right\|_2.
\end{eqnarray}
% Note that, for $\theta_i=\alpha_0$, we have $\frac{\partial}{\partial \theta_i} X_t(\theta)=\frac{1}{1-\beta} \leq \frac{1}{\epsilon}$.
 For this model, the function $\frac{\partial}{\partial \theta}f^{\theta}_\infty$ satisfies the classical Lipschitz-type condition
\[ \left\|\frac{\partial}{\partial \theta_i}f^{\theta}_\infty(y_1,y_2,\cdots)- \frac{\partial}{\partial \theta_i}f^{\theta}_\infty(y^{\prime}_1,y^{\prime}_2,\cdots)\right\|_\Theta \leq \overset{\infty}{\underset{k=1}{\sum}} \alpha^{(1)}_k \left| y_k-y^{\prime}_k\right|\]
for all $(y_1,y_2,\cdots), (y^{\prime}_1,y^{\prime}_2,\cdots) \in (\R^{+})^{\N}$ and $i=1,2,3$ where $(\alpha^{(1)}_k)_{k \geq 1}$ is a non-negative sequence satisfying
$\overset{\infty}{\underset{k=1}{\sum}} \alpha^{(1)}_k < + \infty$ (see for instance \cite{Doukhan_CH_P}).\\
Therefore, by choosing $y^{\prime}_k=0$ for all $k$, we get
\begin{eqnarray*}
\left\| \frac{\partial}{\partial \theta_i} X_t(\theta)\right\|_\Theta = \left\|\frac{\partial}{\partial \theta_i}f^{\theta}_\infty(Y_{t-1},Y_{t-2},\cdots) \right\|_\Theta
\leq \left\|\frac{\partial}{\partial \theta_i}f^{\theta}_\infty(0,0,\cdots) \right\|_\Theta + \overset{\infty}{\underset{k=1}{\sum}} \alpha^{(1)}_k Y_{t-1}
\leq C+ \overset{\infty}{\underset{k=1}{\sum}} \alpha^{(1)}_k Y_{t-1}.
\end{eqnarray*}
Hence, we have
\begin{eqnarray*}
\left\|\left\| \frac{\partial}{\partial \theta_i} X_t(\theta)\right\| _\Theta\right\|_2 &\leq&  \left\|C+ \overset{\infty}{\underset{k=1}{\sum}} \alpha^{(1)}_k Y_{t-1}\right\|_2
\leq C+\overset{\infty}{\underset{k=1}{\sum}} \alpha^{(1)}_k \left\| Y_{t-1}\right\|_2
\leq C+\left(\E \big(Y^2_{0} \big)\right)^{1/2} \overset{\infty}{\underset{k=1}{\sum}} \alpha^{(1)}_k < \infty.
\end{eqnarray*}
Similarly, we have $\left\|\left\|\frac{\partial}{\partial\theta_j}f^{\theta}_\infty\right\|_\Theta\right\|_2  < \infty$.
 Hence, from (\ref{Major1}),
  \[\E \left[ \sup_{\theta \in \Theta} \left| A^{\prime \prime}(\eta_t(\theta)) \left(\frac{\partial}{ \partial \theta_i}\eta_t(\theta)\right) \left(\frac{\partial}{ \partial \theta_j}\eta_t(\theta)\right) \right| \right]<\infty. \]
%Hence $\E \left[ \sup_{\theta \in \Theta} \left| A^{\prime \prime}(\eta_t(\theta))\left(\frac{\partial}{\partial \theta_i}\eta_t(\theta)\right)^{2} \right| \right]<\infty$.\\
By using the same techniques, we get $ \E \big[ \sup_{\theta \in \Theta} | \big(Y_{t}-A^{\prime}(\eta_t(\theta ))\big) (\partial^{2} \eta_t(\theta) / \partial\theta_i \partial\theta_j ) | \big] < + \infty $,
           for any $i,j = 1,\cdots, d$ ; which shows that (\textbf{A6}) holds. \\
%the condition (\textbf{A6}) is satisfied because $A^{\prime \prime}(\eta_t)=e^{\eta_t}=X_t \geq \frac{\alpha_{0}}{1-\beta}$.\\

\item Consider the negative binomial INGARCH$(1, 1)$  (NB-INGARCH$(1, 1)$) model defined by
% The second linear example of the class (\ref{Model}) that we will consider is the negative binomial INGARCH$(1, 1)$  (NB-INGARCH$(1, 1)$) model given by
\begin{eqnarray}\label{NB-INGARCH(1,1)}
    Y_{t}|\mathcal{F}_{t-1}\sim\textrm{NB}(r,p_{t}), \text{ with } r\frac{(1-p_{t})}{p_{t}} = \mathbb{E}(Y_t|\mathcal{F}_{t-1}) = X_t = \alpha^*_{0} +\alpha^* Y_{t-1} + \beta^* X_{t-1} ;
    \end{eqnarray}
  the true parameter $\theta_0 = (\alpha_0^*, \alpha^* , \beta^*)$ belongs to a compact set $\Theta \subset (0, +\infty) \times [0, +\infty)^2$ such that
 $\alpha + \beta < 1$ for all $\theta = (\alpha_0, \alpha, \beta) \in \Theta$, $r \in \mathbb{N}$ and NB$(r,p)$ denotes the negative binomial distribution with parameters $r$ and $p$.
 This model is considered by assuming that $r$ is fixed and supposed to be known.
 It is slightly different to that defined by Zhu \cite{Zhu_2011} where the (auto)regression have been done according to $\frac{1-p_t}{p_t}=\lambda_t$.
 Model (\ref{NB-INGARCH(1,1)}) is a particular case of (\ref{Model}) with $\eta_t=\log\left(\frac{X_t}{X_t+r}\right)$ and $A(\eta_t)=r\log\left(\frac{r}{1-e^{\eta_t}}\right)$.
  It is easily seen that assumption (\textbf{A}($\Theta$)) holds.
  If the compact set $\Theta$ is defined as in (\ref{Set_Theta}), then one can go along similar lines as in the  model (\ref{INGARCH(1,1)}) to show that the conditions (\textbf{A1})-(\textbf{A6}) are satisfied.
 \end{enumerate}

 \subsubsection{Binary time series}
 Let $\left(Y_t\right)_{t\in \mathbb{Z}}$ be a binary ($0,1$ valued) time series satisfying :
    \begin{eqnarray}\label{Linear_Binary(1,1)}
    Y_{t}|\mathcal{F}_{t-1}\sim\textrm{B}(X_{t}) \text{ with } X_{t} = \alpha^*_{0} +\alpha^* Y_{t-1} + \beta^* X_{t-1} ;
    \end{eqnarray}
   the true parameter $\theta_0 = (\alpha_0^*, \alpha^* , \beta^*) \in \Theta$ where $\Theta$ is a compact subset of $(0, +\infty) \times [0, +\infty)^2$ such that
     $\alpha_0 + \alpha + \beta < 1$ for all $\theta = (\alpha_0, \alpha, \beta) \in \Theta$, B$(X_{t})$ denotes the Bernoulli distribution of parameter $X_t$.
    This is a linear example of the class (\ref{Model}), with $\eta_t=\log\left(\frac{X_{t}}{1-X_{t}}\right)$ and $A(\eta_{t})=\log\left(1+e^{\eta_{t}}\right)$.
    In this model, $Y_{t}|\mathcal{F}_{t-1}$ follows a Bernoulli distribution such that the success probability is a function of the immediate past.
    When the parameter $ \beta^* $ is equal to zero, we obtain a particular case of the binary autoregressive models studied by Hudecov{\'a} \cite{Hudecova_2013}.
    Fokianos {\it et al.} \cite{Fokianos_2013a} have also studied a similar model which the conditional mean depends on a vector of a covariates.
    It is easily seen that assumption (\textbf{A}($\Theta$)) holds.
  Define the compact set $\Theta$ by :
 \[
 \Theta=\Big\{\theta=(\alpha_{0},\alpha,\beta) \in \mathbb{R}^{3}_{+}: ~\epsilon \leq \alpha_{0} + \alpha + \beta \leq 1- \epsilon \Big\}~\rm{for~ some}~ \epsilon>0.
 \]
 Then, by going along similar lines as in the model (\ref{INGARCH(1,1)}), the assumptions (\textbf{A1})-(\textbf{A6}) holds.

\subsubsection{Threshold autoregressive models}
\begin{enumerate}
\item Firstly, we consider a threshold Poisson autoregressive model (INTARCH$(1,1)$) defined by (see also Doukhan and Kengne (2015)):
\begin{eqnarray}\label{INTARCH(1,1)}
    Y_{t}|\mathcal{F}_{t-1}\sim\textrm{Poisson}(X_{t}), \text{ with } X_{t}=\alpha^*_{0} +\alpha^*_{1}\max (Y_{t-1}-\ell,0)+\alpha^*_{2} \min (Y_{t-1},\ell) + \beta^* X_{t-1} ;
    \end{eqnarray}
  the true parameter $\theta_0 = (\alpha_0^*, \alpha^*_1, \alpha^*_2 , \beta^*) \in \Theta$ where $\Theta$ is a compact subset of $(0, +\infty) \times [0, +\infty)^3$
     such that $ \max(\alpha_1, \alpha_2) + \beta < 1$ for all $\theta=(\alpha_{0},\alpha_1, \alpha_2 ,\beta) \in \Theta$
  and $\ell$ is a non-negative integer valued, called the threshold parameter of the model.
%This nonlinear model is obtained by using the definition like the threshold ARCH model proposed by Zakoian (1994)  for to reflect the asymmetric relation between volatility and past returns.
%The parameter $l$ is called the threshold of the model.
 The condition (\textbf{A}($\Theta$)) is satisfied.  We have the following representation for all $\theta \in \Theta $ :
 \[ X_{t}(\theta)=\frac{\alpha_{0}}{1-\beta} + \sum_{k=0}^{\infty}\beta^{k} \Big(\alpha_1 \max (Y_{t-1-k}-\ell,0)+ \alpha_2 \min (Y_{t-1-k},\ell)\Big):=f^{\theta}_\infty(Y_{t-1},\cdots). \]
 %where $\left(\beta + \max (\alpha_{1},\alpha_{2})\right)<1$.
 Hence, $f^{\theta}_\infty(Y_{t-1},\cdots) \geq \frac{\alpha_{0}}{1-\beta} := x^{*}_{\theta}$.
 Define the compact set $\Theta$ by :
 \[
\Theta=\Big\{\theta=(\alpha_{0},\alpha_1, \alpha_2 ,\beta) \in \mathbb{R}^{4}_{+}:~0<\alpha_{L}\leq \alpha_{0} \leq \alpha_{U},~\epsilon \leq \max(\alpha_1, \alpha_2) + \beta \leq 1- \epsilon \Big\}~\rm{for~ some}~ \alpha_L, \alpha_U, \epsilon>0.
\]
 Then, by going along similar lines as in the model (\ref{INGARCH(1,1)}), the assumptions (\textbf{A1})-(\textbf{A6}) holds.

\item [$2.$] We assume that the process $(Y_t)_{t \in \mathbb{Z}}$ follows a negative binomial INTARCH model (NB-INTARCH) defined by
  \begin{align}
   \nonumber  Y_{t}|\mathcal{F}_{t-1} & \sim \text{NB}(r,p_{t}), \text{ with }  \\
  \label{NB-INTARCH(1,1)}   & r\frac{(1-p_{t})}{p_{t}} = \mathbb{E}(Y_t|\mathcal{F}_{t-1}) = X_{t}=\alpha^*_{0} +\alpha^*_{1}\max (Y_{t-1}-\ell,0)+\alpha^*_{2} \min (Y_{t-1},\ell) + \beta^* X_{t-1} ;
    \end{align}
  the true parameter $\theta_0 = (\alpha_0^*, \alpha^*_1, \alpha^*_2 , \beta^*) \in \Theta$ where $\Theta$ is a compact subset of $(0, +\infty) \times [0, +\infty)^3$
     such that $ \max(\alpha_1, \alpha_2) + \beta < 1$ for all $\theta=(\alpha_{0},\alpha_1, \alpha_2 ,\beta) \in \Theta$, $\ell$ the threshold parameter of the model
    and NB$(r,p)$ denotes the negative binomial distribution of parameter $(r,p)$ where $r$ is fixed and supposed to be known.
  The assumption (\textbf{A}($\Theta$)) holds.
   Define $\Theta$ as above and one can show that the conditions (\textbf{A1})-(\textbf{A6}) are satisfied.
\end{enumerate}

\section{Change-point test and asymptotic results}
We propose a change-point test based on the maximum likelihood estimator of the parameter of the model (\ref{Model}).
 Using the idea of Doukhan and Kengne \cite{Doukhan_CH_P} and Kengne \cite{Kengne_2012}, we will construct a test statistic that converges to a known distribution under $\textrm{H}_0$
 and diverges to infinity under the alternative of change in the model. Throughout the sequel, the following notations will be used:
\begin{enumerate}
\item $ \|x \|:= \sqrt{\sum_{i=1}^{p} |x_i|^2 } $ for any $x \in \mathbb{R}^{p}$;
\item  $\left\|f\right\|_{\Theta}:=\sup_{\theta \in \Theta}\left(\left\|f(\theta)\right\|\right)$ for any function $f:\Theta \longrightarrow \mathbb{R}^{d^{\prime}}$;
\item  $T_{k,k^{\prime}}:=\left\{k,k+1,\cdots,k^{\prime}\right\}$ for any $k,k^{\prime} \in \left\{1,2,\cdots,n\right\}$ such as $k \leq k^{\prime}$;
\item  $L_n({T_{k,k^{\prime}}},\theta):=\sum_{t=k}^{k^{\prime}}\ell_t(\theta)$ is the conditional log-likelihood function computed on the observations
$Y_k,Y_{k+1},\cdots,Y_{k^{\prime}}$ where  $\ell_t(\theta)=\{\eta_t(\theta )Y_{t}-A^{\prime}(\eta_t(\theta ))\}$ with $\eta_t(\theta)=(A^{\prime})^{-1}(X_t(\theta))$;
\item $\widehat{\theta}_n({T_{k,k^{\prime}}})=\textrm{argmax}_{\theta \in \Theta}\left(L_n({T_{k,k^{\prime}}},\theta)\right)$ is the MLE computed on the observations $Y_k,Y_{k+1},\cdots,Y_{k^{\prime}}$.
\end{enumerate}
Let us recall that, under H$_0$, the asymptotic covariance matrix of $\widehat{\theta}_n({T_{1,n}})$
is
 \[\widehat{\Omega}_{n}=\left.\frac{1}{n}\sum_{t=1}^{n}\left({A}^{\prime \prime}\left(\eta_t(\theta)\right)\left(\frac{\partial \eta_t(\theta)}{\partial \theta}\right)\left(\frac{\partial \eta_t(\theta)}{\partial \theta}\right)^{T}\right)\right|_{\theta=\widehat{\theta}_n({T_{1,n}})}\]
and that it is a consistent estimator of $\Omega$ (see (\ref{def_Sigma})).
But this consistency is not ensured in general under the change-point alternative (see \cite{Doukhan_CH_P}). \\
Let $(u_n)_{n\geq 1}$ and $(v_n)_{n\geq 1}$ be two integer valued sequences satisfying ~$u_n,v_n \rightarrow +\infty$ and $\frac{u_n}{n},\frac{v_n}{n} \rightarrow 0$ as $n\rightarrow +\infty$.
For all $n \geq 1$, define the matrix $\widehat{\Omega}_{n}(u_n)$  by
\begin{multline*}
\widehat{\Omega}_{n}(u_n)=\frac{1}{2}\left[\frac{1}{u_n}\sum_{t=1}^{u_n} {A}^{\prime \prime}\left(\eta_t(\theta)\right)\left(\frac{\partial \eta_t(\theta)}{\partial \theta}\right)\left(\frac{\partial \eta_t(\theta)}{\partial \theta}\right)^{T} \Big|_{\theta=\widehat{\theta}_n(T_{1,u_n})} \right. \\
   \left. + \frac{1}{n-u_n}\sum_{t=u_n+1}^{n} {A}^{\prime \prime}\left( \eta_t(\theta)\right)\left(\frac{\partial \eta_t(\theta)}{\partial \theta}\right)\left(\frac{\partial \eta_t(\theta)}{\partial \theta}\right)^{T} \Big|_{\theta= \widehat{\theta}_n(T_{u_n+1,n})} \right] .
 \end{multline*}
 Under H$_0$, according to the consistency of $\widehat{\Omega}_{n}$ and the ergodicity of the process $(X_t, Y_t)_{t \in \Z}$, $\widehat{\Omega}_{n}(u_n)$ is also a consistent estimator
 of $\Omega$. See \cite{Doukhan_CH_P}  for the motivation of the use of $\widehat{\Omega}_{n}(u_n)$ instead of $\widehat{\Omega}_{n}$.
 Then, consider the test statistic :
 \[ \widehat{C}_n=\max_{v_n\leq k \leq n-v_n}\widehat{C}_{n,k} \]
with
 \[ \widehat{C}_{n,k}=\frac{1}{q^{2}(\frac{k}{n})}\frac{k^{2}(n-k)^{2}}{n^{3}}\left(\widehat{\theta}(T_{1,k})-\widehat{\theta}(T_{k+1,n})\right)^{T}\widehat{\Omega}_{n}(u_n)\left(\widehat{\theta}(T_{1,k})-\widehat{\theta}(T_{k+1,n})\right) \]
 where $q:(0, 1)\rightarrow (0, \infty)$ is a weight function non-decreasing in a neighborhood of zero, non-increasing in a neighborhood of one and satisfying
  \[ \inf_{\varphi<\tau<1-\varphi}q(\tau)>0 ~~\textrm{for all}~~ 0<\varphi< \frac{1}{2}. \]
Its  behavior can be controlled at the neighborhood of zero and one by the integral
\[I(q,c)=\int_{0}^{1}\frac{1}{t(1-t)}\exp\left(-\frac{cq^{2}(t)}{t(1-t)}\right)dt ,~~c>0 \]
see \cite{Csorgo_1986}. The weight function $q$ is used to increase the power of the test based on the statistic $\widehat{C}_n$. \\
The following two theorems give the asymptotic behavior of the statistic $\widehat{C}_n$ under the null and alternative hypothesis.
\begin{thm}\label{Thm_Asy_H0}
Assume that (\textbf{A0})-(\textbf{A6}) and (\textbf{A}($\Theta$)) hold. Under H$_0$,  if there exists $ c>0$ such that $I(q,c)<\infty$, then
\[ \widehat{C}_n \limiteloin \sup_{0\leq\tau\leq1}\frac{\left\|W_d(\tau)\right\|^{2}}{q^{2}(\tau)}  \]
where $W_d$ is a $d$-dimensional Brownian bridge.
\end{thm}
Then, at a nominal level $\alpha \in (0,1)$, the critical region of the test is $(\widehat{C}_{n}>c_\alpha)$, where $c_\alpha$ is the $(1-\alpha)$-quantile of the distribution
of $\sup_{0<\tau<1} \big( \left\|W_d(\tau)\right\|^{2} / q^{2}(\tau) \big)$. Hence, the proposed test has correct size asymptotically.
 When $q\equiv 1$, the values of $c_\alpha$ can be obtained in Lee \textit{et al.} (2003) for $d \in \left\{1,2,\cdots,10\right\}$.

\medskip

 Under the alternative, we assume

  \medskip

\noindent {\bf Assumption B }: {\em there exist $\tau \in (0,1)$ such that $t^{*}=[n\tau^*]$ (where $[x]$ is the integer part of $x$).}

\begin{thm}\label{Thm_Asy_H1}
Under H$_1$, assume that B, (\textbf{A0})-(\textbf{A6}) and (\textbf{A}($\Theta$)) hold. Then,
\[ \widehat{C}_n  \limiteproban +\infty .\]
\end{thm}
 Therefore, the proposed procedure based on $\widehat{C}_n$ is consistent in power.
 Under H$_1$, a classical estimator of the breakpoint is
 \[   \widehat{t}_n =  \underset{v_n \leq k  \leq n - v_n }{ \text{argmax}  }  \widehat{C}_{n,k} . \]

\section{Some simulations results}
In this section, we evaluate the performance of the proposed procedure through an empirical study. We consider a linear and a nonlinear example of the class of models
(\ref{Model}). For a sample size $n$, the test statistic $\widehat{C}_n$ is computed with $q\equiv 1$ and $u_n=v_n=[\left(\log(n)\right)^{2}]$.
The nominal level considered in the sequel is $\alpha = 0.05$. We denote by $\theta = (\alpha_{0},\alpha,\beta)$ the parameter of the model considered.
\subsection{Test for parameter change in NB-INGARCH$(1,1)$ models}
\begin{figure}[h!]
\begin{center}
\includegraphics[height=11.85cm, width=16.75cm]{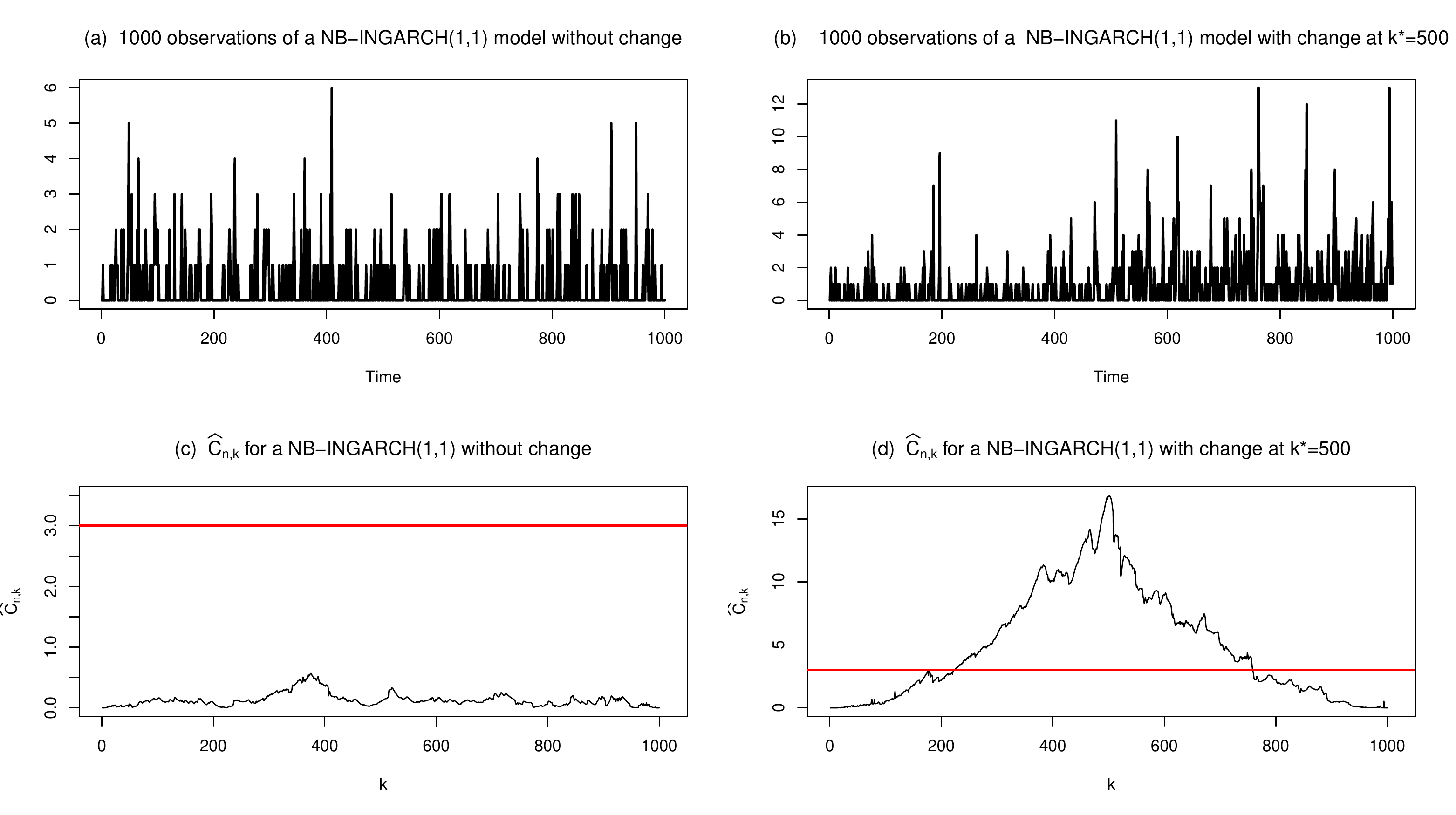}
\end{center}
\caption{Typical realization of $1000$ observations of two NB-INGARCH$(1,1)$ processes and the corresponding
statistics $\widehat{C}_{n,k}$. \textbf{(a)} is a NB-INGARCH$(1,1)$ process without change, where the parameter
$\theta_0=(0.20, 0.30, 0.25)$ is constant. \textbf{(b)} is a NB-INGARCH$(1,1)$ process where the parameter $\theta_0=(0.20, 0.30, 0.25)$ changes
to $\theta_1=(0.70, 0.3, 0.25)$ at $k^{*}=500$. \textbf{(c)} and \textbf{(d)} represent their corresponding statistics $\widehat{C}_{n,k}$.
The horizontal line represents the limit of the critical region of the test.}
\label{Graphe_Sim_NB_INGARCH}
\end{figure}

Consider the negative binomial INGARCH$(1,1)$ model given by
\begin{eqnarray}\label{Model_Bin_INGARCH}
    Y_{t}|\mathcal{F}_{t-1}\sim\textrm{NB}(r,p_{t}), \text{ with }   r\frac{(1-p_{t})}{p_{t}}=X_{t}=\alpha^*_{0} +\alpha^* Y_{t-1} + \beta^* X_{t-1}
    \end{eqnarray}
 where $\alpha^*_{0}>0$, $\alpha^*, \beta^* \geq 0$, see (\ref{NB-INGARCH(1,1)}).
 For the problem of estimating the parameter $r$, we can use an information criteria such as AIC or BIC (see Davis and Wu (2009)  or Zhu (2011)).\\
Now, assume that the parameter $r$ is known and consider the problem of testing for parameter change in the model (\ref{Model_Bin_INGARCH}).
For $r=1$ and $n=1000$, we consider two trajectories $(Y_1, \cdots, Y_n)$ generated from (\ref{Model_Bin_INGARCH})
in the following situations : a scenario without change and a scenario with a change at $k^{*} = 500$.
Figure \ref{Graphe_Sim_NB_INGARCH} indicates the statistic $\widehat{C}_{n,k}$ of the test for these two scenarios.
 One can see that, under the null hypothesis (i.e. no change), the statistic $\widehat{C}_{n}$  is less than the limit of the critical region which is represented by
 the horizontal line (see Figure \ref{Graphe_Sim_NB_INGARCH}\textbf{(a)} and Figure \ref{Graphe_Sim_NB_INGARCH}\textbf{(c)}).
 Under the alternative, the statistic $\widehat{C}_{n,k}$ is large around the point where the change occurs and $\widehat{C}_{n}$ is greater
 than the critical value of the test (see Figure \ref{Graphe_Sim_NB_INGARCH}\textbf{(b)} and Figure \ref{Graphe_Sim_NB_INGARCH}\textbf{(d)}).

For a sample size $n = 500,1000$, Table \ref{Levels_powers_BINOMIAL_INGARCH} indicates the empirical levels computed when the parameter
is $\theta_0$ and the empirical powers computed when $\theta_0$ changes to $\theta_1$ at half the sample size, where the number of failures $r=1,8$ is fixed ;
these results are obtained after $200$ replications.
The choice of $r=8$ is related to the real data example (see below).
These results indicate that the procedure produces a  reasonable empirical levels which
approaching the nominal one when $n$ increases, but these results are more accurate for $r=1$.
The empirical powers of the procedure are also well satisfactory.

\begin{table}[h!]
\centering
\begin{tabular}{|cc|c|c|c|}
\hline
&&&$n=500$&$n=1000$\\
\hline
\hline
Empirical levels : & & &&\\
                  &$\theta_0=(1.0,0.2,0.15)$ &$r=1$&$0.030$&$0.025$\\\
                  &&$r=8$&$0.015$&$0.020$\\
& & &&\\
                  & $\theta_0=(0.2,0.3,0.25)$ &$r=1$&$0.075$&$0.045$\\
                  &&$r=8$&$0.040$&$0.020$\\
& & &&\\
                  & $\theta_0=(8, 0.2, 0.15)$ &$r=1$&$0.035$&$0.040$\\
                  &&$r=8$&$0.025$&$0.030$\\\
&&&&\\
\hline
Empirical powers : & & &&\\
                  &$\theta_0=(1.0,0.2,0.15)~; ~\theta_1=(1.0,0.05,0.85)$ &$r=1$&$0.985$&$0.995$\\
                  &&$r=8$&$0.995$&$0.995$\\
& & &&\\
                  & $\theta_0=(0.2,0.3,0.25)~; ~\theta_1=(0.7,0.3,0.25)$ &$r=1$&$0.980$&$0.990$\\
                  &&$r=8$&$0.990$&$0.995$\\
& & &&\\
                  & $\theta_0=(8, 0.2, 0.15)~; ~\theta_1=(2, 0.15, 0.5)$ &$r=1$&$0.945$&$0.980$\\
                  &&$r=8$&$0.985$&$0.995$\\
 \hline
\end{tabular}
\caption{{\footnotesize Empirical levels and powers at the nominal level $0.05$ of test for parameter change in a NB-INGARCH$(1,1)$ model
with one change-point alternative.}}
\label{Levels_powers_BINOMIAL_INGARCH}
%\scriptsize
\end{table}

\subsection{Test for parameter change in binary time series models}
We consider a binary ($0,1$ valued) process $\left(Y_t\right)_{t\in \mathbb{Z}}$ satisfying :
    \begin{eqnarray}\label{Model_Binary(1,1)}
    Y_{t}|\mathcal{F}_{t-1}\sim\textrm{B}(X_{t})~~\textrm{with}~~X_{t}=\alpha^*_{0} +\alpha^* Y_{t-1} + \beta^*  X_{t-1}
    \end{eqnarray}
 with  $\alpha^*_{0} > 0$; $\alpha^*, \beta^* \geq 0$, see (\ref{Linear_Binary(1,1)}).
 We consider the problem of testing for parameter change in this model.
 Let us recall that Hudecov{\'a} \cite{Hudecova_2013} has studied the problem of testing for
 parameter change in a similar model. They have used the logit link function to take into account the dependence between the success probability and the previous values of
 the series. They proposed a change-point test based on a normalized cumulative sums of residuals. We have made a comparison with their procedure
 on the real data application (see bellow, the application to the US recession data).

\begin{figure}[h!]
\begin{center}
\includegraphics[height=11.95cm, width=17cm]{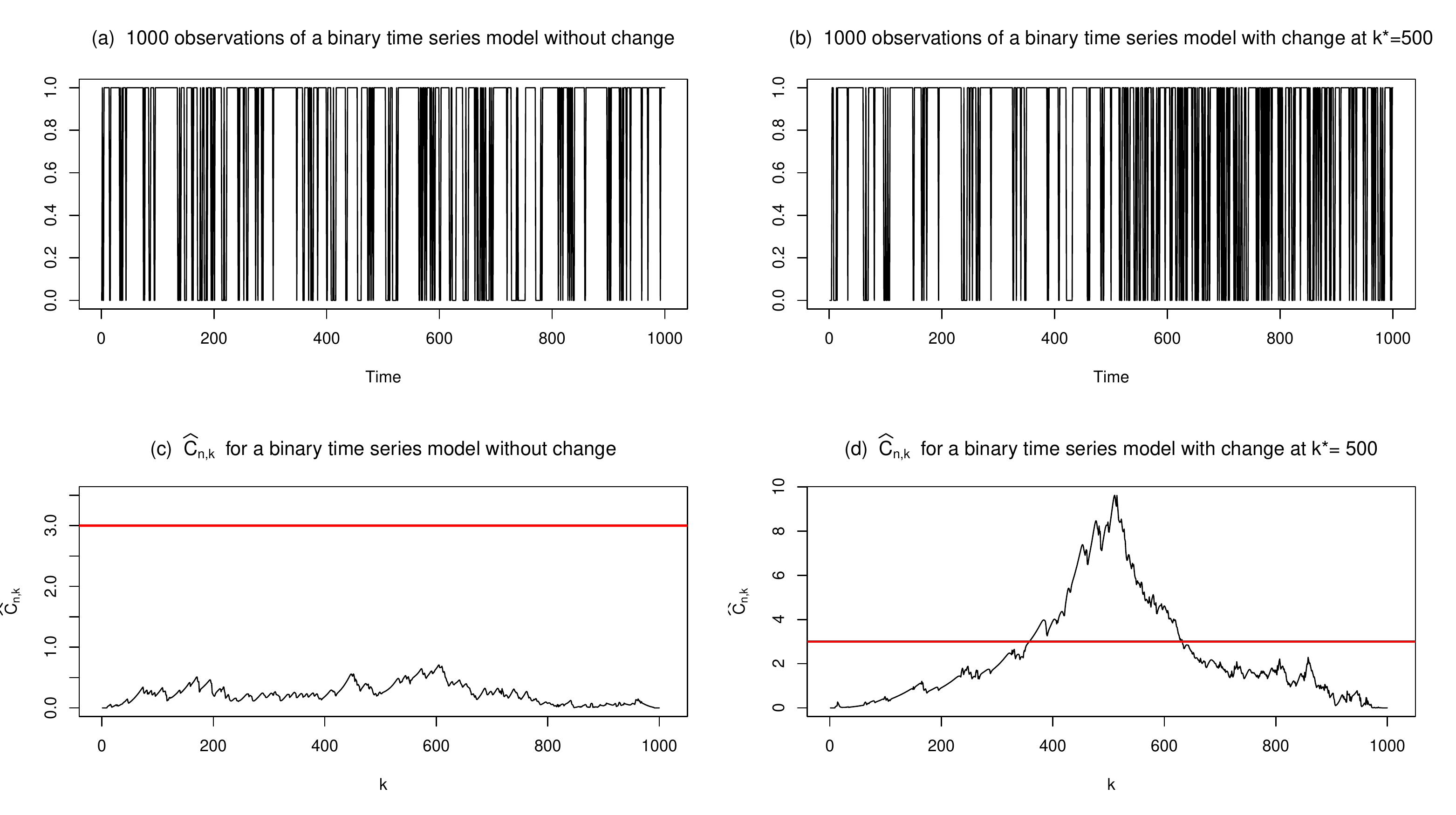}\vspace{-0.30cm}
\end{center}
\caption{Typical realization of $1000$ observations of two binary time series processes and the corresponding
statistics $\widehat{C}_{n,k}$. \textbf{(a)} is a binary time series process without change, where the parameter
$\theta_0=(0.2, 0.35, 0.4)$ is constant. \textbf{(b)} is a binary time series process where the parameter $\theta_0=(0.2, 0.35, 0.4)$ changes
to $\theta_1=(0.35, 0.1, 0.4)$ at $k^{*}=500$. \textbf{(c)} and \textbf{(d)} represent their corresponding statistics $\widehat{C}_{n,k}$.
 The horizontal line represents the limit of the critical region of the test.}
 \label{Graphe_Sim_BINARY_INGARCH}
\end{figure}

For $n=1000$, Figure \ref{Graphe_Sim_BINARY_INGARCH} indicates a typical realization of the statistic $\widehat{C}_{n,k}$ for a series generated
from (\ref{Model_Binary(1,1)}). We consider a scenario without change (see Figure \ref{Graphe_Sim_BINARY_INGARCH}\textbf{(a)}) and a scenario with a change at $k^{*} = 500$
(see Figure \ref{Graphe_Sim_BINARY_INGARCH}\textbf{(b)}).
 One can see that, in a scenario without change, the test statistic $\widehat{C}_{n}$ is under the horizontal line which represents the limit of the critical region
 (see Figure \ref{Graphe_Sim_BINARY_INGARCH}\textbf{(c)}). Under the alternative (one change),
 this statistic is greater than the critical value in the neighborhood of the breakpoint $k^{*} = 500$ (see Figure \ref{Graphe_Sim_BINARY_INGARCH}\textbf{(d)}).\\

For $n =500,1000$; we generate a trajectory $(Y_1,Y_2,\cdots,Y_n)$ in the following situations :  the parameter $\theta=\theta_0$ is constant (no change) and
 the parameter changes from $\theta_0$ to $\theta_1$ at $n/2$ (one change). For different values of  $\theta_0$ and $\theta_1$, Table \ref{Levels_powers_BINARY_INGARCH}
indicates the empirical levels  and the empirical powers based of $200$ replications.
 Again, these results display the accuracy of the proposed procedure.

\begin{table}[h!]
\centering
\begin{tabular}{|ccc|c|c|}
\hline
&&&$n=500$&$n=1000$\\
\hline
\hline
Empirical levels :&&&&\\

                  &$\theta_0=(0.1,0.75,0.05)$&&$0.045$&$0.050$\\
& & & & \\
                  &$\theta_0=(0.2, 0.35, 0.4)$&&$0.035$&$0.045$\\
&&& & \\
\hline
Empirical powers :&&&&\\
                  &$\theta_0=(0.1,0.75,0.05)$ & $\theta_1=(0.05,0.7,0.05)$&$0.605$&$0.820$\\
                  & & $\theta_1=(0.1,0.4,0.05)$&$0.995$&$0.995$\\
& & && \\
                  &$\theta_0=(0.2, 0.35, 0.4)$ & $\theta_1=(0.35, 0.1, 0.4) $&$0.775$&$0.960$\\ %%S2
                  & & $\theta_1=(0.05, 0.35, 0.4) $&$0.985$&$0.990$\\
 \hline
\end{tabular}
\caption{{\footnotesize Empirical levels and powers at the nominal level $0.05$ of test for parameter change in binary time series models
with one change-point alternative.}}
\label{Levels_powers_BINARY_INGARCH}
\end{table}

\section{Real data application}
\subsection{Number of transactions of Ericsson B stock}
\begin{figure}[h!]
\begin{center}
\includegraphics[height=7cm, width=13cm]{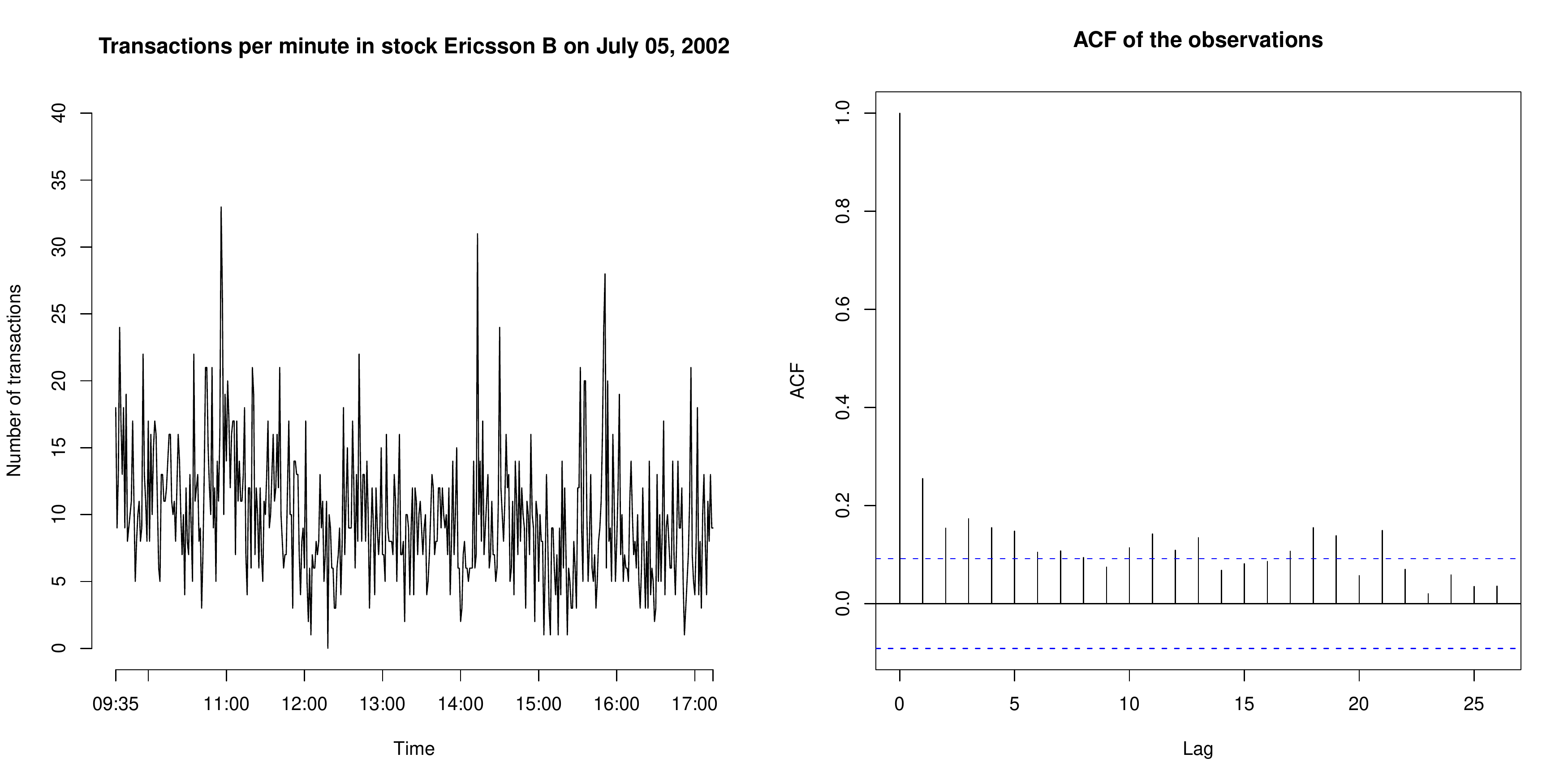}\vspace{-0.30cm}
\end{center}
\caption{{\footnotesize Number of transactions per minute for the stock Ericsson B during July 05, 2002 and their autocorrelation
function.}}
 \label{Plot_Ericsson_05}
\end{figure}

We apply this change-point test to the series of the number of transactions per minute for the stock Ericsson B during July 05, 2002.
There are $460$ available observations which represent the transaction of approximately $8$ hours (from $09:35$ through $17: 14$).
Fokianos \textit{et al.} \cite{Fokianos_2009} and Davis and Liu \cite{Davis_2012} have fitted the series of July 2, 2002 by using the Poisson INGARCH$(1,1)$ and the negative binomial INGARCH$(1,1)$ model.
The plots of the data (the series of July 05, 2002) and the autocorrelation function displayed on Figure \ref{Plot_Ericsson_05} show a positive dependence.
The empirical mean and variance of this series are $9.824$ and $23.753$ respectively. Hence this data are overdispersed. The positive dependence and the overdispersion observed
suggest that the NB-INGARCH model is a candidate to fit these data.\\
We consider the model (\ref{Model_Bin_INGARCH}) with $r=8$, and we test whether a change has occurred in the parameter in the model for these data.
The value $r=8$ is that have been obtained by Davis and Liu \cite{Davis_2012} on the series of July 2, 2002.
The test statistic is computed with $q \equiv 1$ and $u_n=v_n=[\left(\log(n)\right)^{2}]$. We set the initial value $X_1= 9.824$ (the empirical mean of the series).

\begin{figure}[h!]
\begin{center}
\includegraphics[height=10.5cm, width=17.5cm]{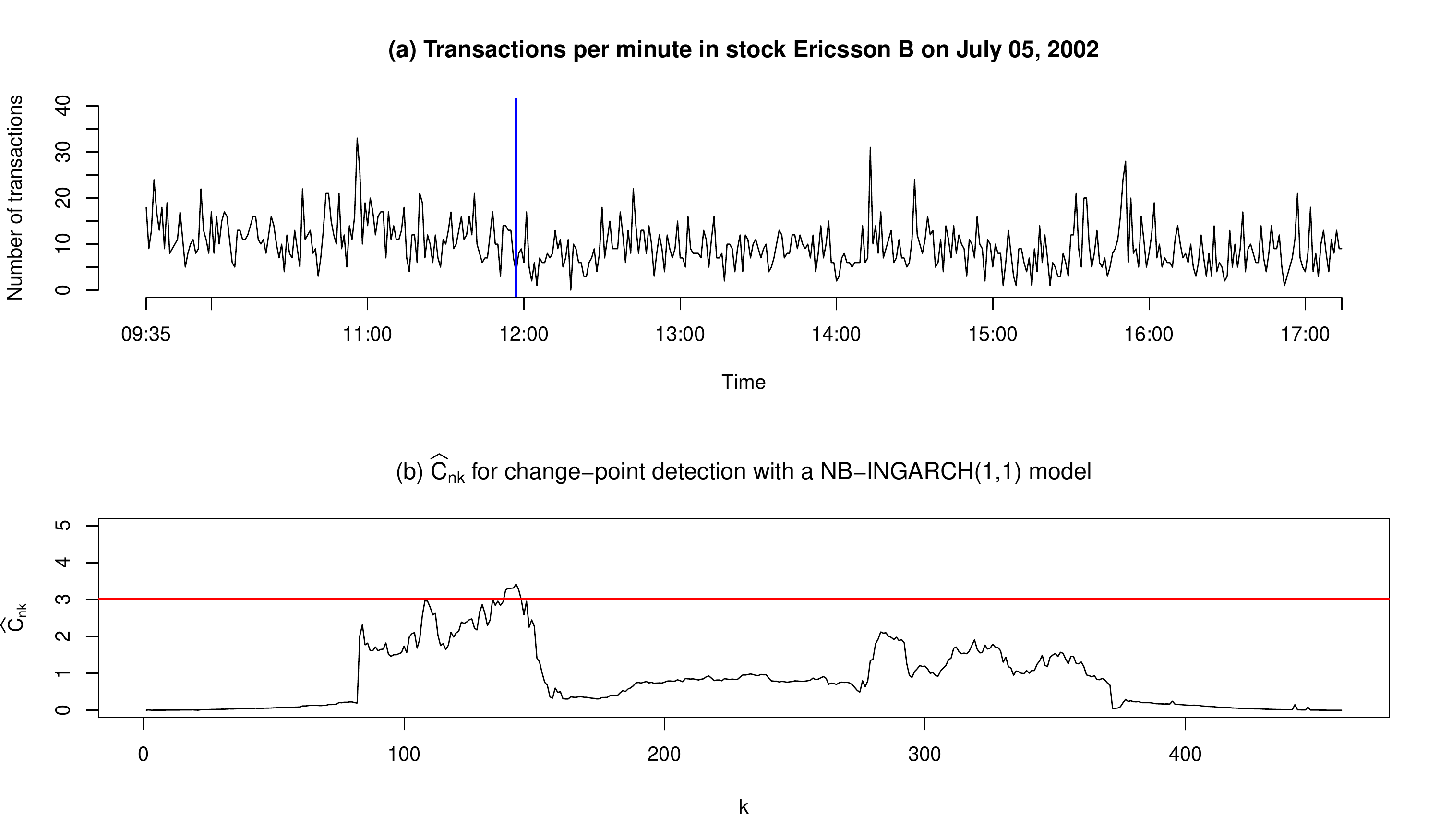}
\end{center}
\caption{{\footnotesize Plot of $\widehat{C}_{n,k}$ of the number of transactions per minute for the stock Ericsson B during July 05, 2002 with NB-INGARCH$(1,1)$ models,
where the horizontal line represents the limit of the critical region of the test and the vertical line is the estimated breakpoint.}}
 \label{Plot_Stat_Ericsson_05}
\end{figure}

The corresponding critical value on nominal level $\alpha=0.05$ is $c_\alpha=3.004$ and the value of the test statistic $\widehat{C}_n$ is equal to $3.406$.
This indicates a change in the parameter of the model. $\widehat{C}_{n,k}$ has a maximum at $\widehat{t}=143$ (see Figure \ref{Plot_Stat_Ericsson_05}),
that is, the location of the change is $11:57$. The estimated model with a change is given by

\begin{equation*}
X_{t}=\left\{
\begin{array}{ll}
 8.1462 +0.1945  Y_{t-1} + 0.1318 X_{t-1}~\textrm{for}~t\leq 143,\\
\\
2.5902 +0.1253  Y_{t-1} + 0.5794  X_{t-1}~\textrm{for}~t > 143.
\end{array}
\right.
\end{equation*}
Such change around the midday have also been obtained by Doukhan and Kengne \cite{Doukhan_CH_P} for the number of transactions during July 16, 2002.

\subsection{US recession data}
\begin{figure}
\begin{center}
\includegraphics[height=9cm, width=14cm]{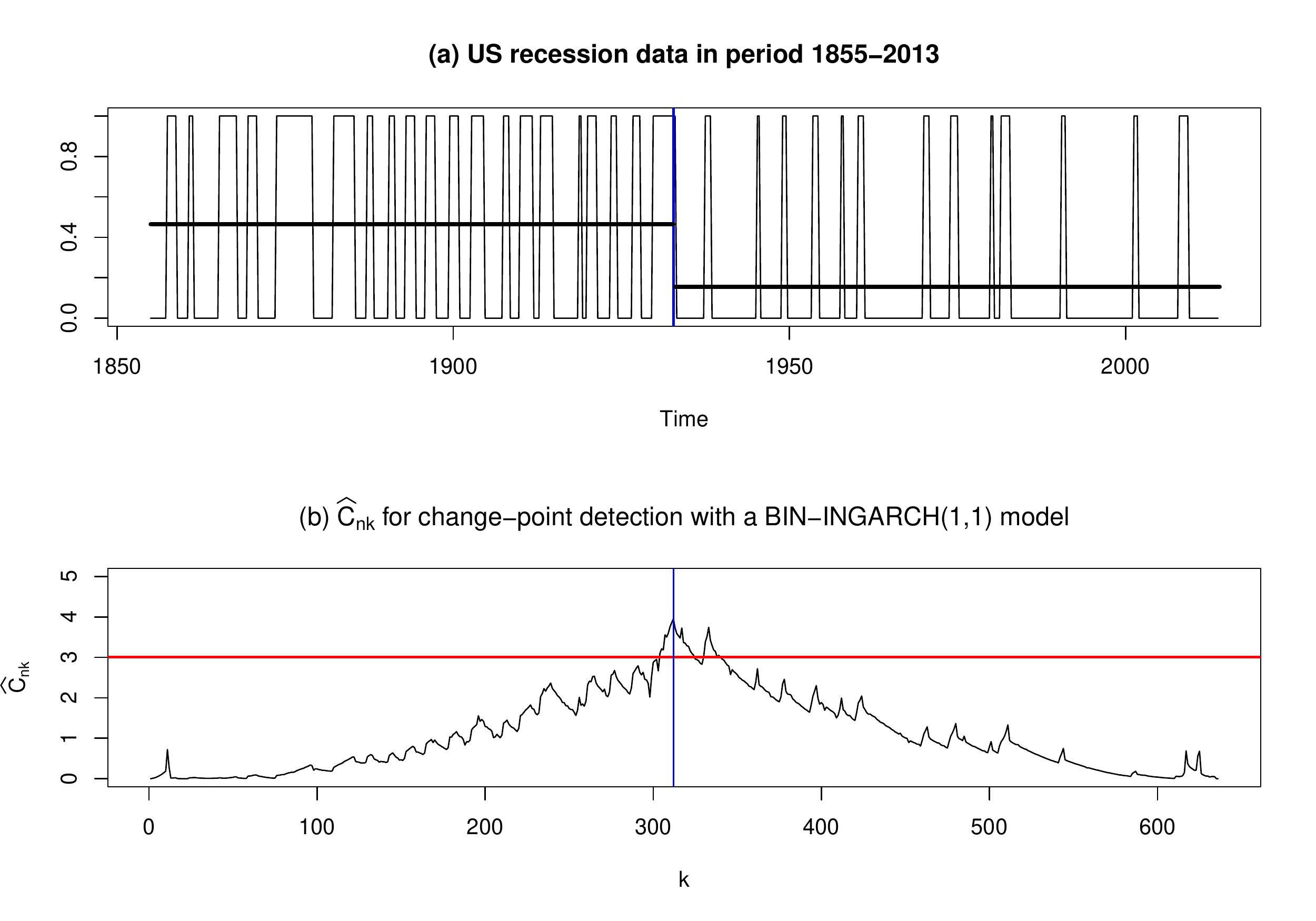} %\vspace{-0.30cm}
\end{center}
\caption{{\footnotesize Plot of $\widehat{C}_{n,k}$ of the US recession data in period 1855-2013 with a binary time series model, where the horizontal line in \textbf{(b)} represents the limit of the critical region of the test,
the vertical line is the estimated breakpoint and the horizontal line in \textbf{(a)} corresponds
to the probabilities of a recession before and after the  breakpoint.}}
\label{Plot_US_Recession}
\end{figure}

As a second example, we consider the quarterly recession data (see Figure \ref{Plot_US_Recession}) from the USA for the period 1855-2013.
In economics, a recession is a slowing phenomenon of the pace of the economic growth reflected in a downswing in the GDP (Gross Domestic Product).
For this series, $Y_t$ is a binary variable that is equal to 1 if there is a recession in at least one month in the quarter, and 0 otherwise.
There are $636$ quarterly observations obtained from The National Bureau of Economic Research. The data are available online at the web site
"http://research.stlouisfed.org/fred2/series/USREC/downloaddata".
These data have been already analyzed by different authors.
Recently, Hudecov{\'a} \cite{Hudecova_2013} has applied a change-point procedure based on a normalized cumulative sums of residuals for these data and has found a break in the first quarter of 1933. \\
Now, we consider the model (\ref{Model_Binary(1,1)}) and apply our procedure to these data.
Figure \ref{Plot_US_Recession} displays the realizations of  $\widehat{C}_{n,k}$ computed with $q \equiv 1$ and $u_n=v_n=[\left(\log(n)\right)^{2}]$.
As shown in this figure, a change has been detected at $ \widehat{t}=312$, which corresponds to the last quarter of 1932. These results are in accordance with those obtained
by Hudecov{\'a} \cite{Hudecova_2013}. This period corresponds to the election of the President Roosevelt in US and the beginning of the end of the great depression.
The estimated model with a change is given by

\begin{equation*}
X_{t}=\left\{
\begin{array}{ll}
0.1193 +0.7483  Y_{t-1} + 1.3614 \times 10^{-8} X_{t-1}~\textrm{for}~t \leq 312,\\
\\
0.0474 +0.6668  Y_{t-1} + 4.0029 \times  10^{-10}  X_{t-1}~\textrm{for}~t> 312.
\end{array}
\right.
\end{equation*}

\section{Proof of the main results}
\subsection{Proof of Theorem \ref{Thm_Asy_H0}}
Define the statistic
$$C_n=\max_{v_n<k<n-v_n}C_{n,k}~~~~~~~~~~~~~~~~~~~~~~~~~~~~~~~~~~~~~~~~~~~~~~~~~~~~~~~~~~$$
where $$C_{n,k}=\frac{1}{q^{2}(\frac{k}{n})}\frac{k^{2}(n-k)^{2}}{n^{3}}\left(\widehat{\theta}(T_{1,k})-\widehat{\theta}(T_{k+1,n})\right)^{T}\Omega \left(\widehat{\theta}(T_{1,k})-\widehat{\theta}(T_{k+1,n})\right).$$~
Let $k,k^{\prime} \in[1,n]$, $\bar{\theta} \in \Theta$ and $i \in \{1,2,\cdots,n\}$. The Taylor expansion to the function $\theta \mapsto \frac{\partial}{\partial \theta_i} L_n({T_{k,k^{\prime}}},\theta)$ implies that there exists $\theta_{n,i}$ between $\bar{\theta}$ and $\theta_{0}$ such that
$$ \frac{\partial}{\partial \theta_i} L_n({T_{k,k^{\prime}}},\bar{\theta})=\frac{\partial}{\partial \theta_i} L_n({T_{k,k^{\prime}}},\theta_0) +\frac{\partial^{2}}{\partial \theta\partial \theta_i} L_n({T_{k,k^{\prime}}},\theta_{n,i})(\bar{\theta}-\theta_0).$$
It is equivalent to
\begin{eqnarray}\label{Eq_Taylor}
(k^{\prime}-k+1)G_n({T_{k,k^{\prime}}},\bar{\theta}).(\bar{\theta}-\theta_0)=\frac{\partial}{\partial \theta} L_n({T_{k,k^{\prime}}},\theta_0)-\frac{\partial}{\partial \theta} L_n({T_{k,k^{\prime}}},\bar{\theta})
\end{eqnarray}
where
\begin{eqnarray}\label{G_n}
G_n({T_{k,k^{\prime}}},\bar{\theta})= - \frac{1}{(k^{\prime}-k+1)}\frac{\partial^{2}}{\partial \theta\partial \theta_i} L_n({T_{k,k^{\prime}}},\theta_{n,i})_{1\leq i \leq d}.
\end{eqnarray}
 The following lemma will be useful in the sequel.
\begin{lem}\label{Lem_H0}
Suppose that the assumptions of the Theorem \ref{Thm_Asy_H0} hold. Then,
\begin{enumerate}
    \item  $\max_{v_n<k<n-v_n}\left|\widehat{C}_{n,k}-C_{n,k}\right|=o_P(1)$;
    \item $\left(\frac{\partial}{\partial \theta} \ell_t(\theta_0),\mathcal{F}_{t}\right)_{t \in \mathbb{Z}}$ is a stationary ergodic, square integrable martingale difference sequence with covariance matrix $\Omega$;
\item $\mathbb{E}\left(\frac{\partial^{2} \ell_0 (\theta_0)}{\partial \theta\partial \theta^{T}}\right)=-\Omega$;
\item $G_n({T_{1,n}},\widehat{\theta}_n({T_{1,n}})) \stackrel{a.s}{\longrightarrow}\Omega$
 and ~
  $ \left.\dfrac{1}{n}\sum_{t=1}^{n}\left({A}^{\prime \prime}\left( \eta_t(\theta)\right)\left(\frac{\partial \eta_t(\theta)}{\partial \theta}\right)\left(\frac{\partial \eta_t(\theta)}{\partial \theta}\right)^{T}\right)\right|_{\theta=\widehat{\theta}_n }
     \stackrel{a.s}{\longrightarrow}\Omega$ as $n\rightarrow +\infty$.

\end{enumerate}
\end{lem}
\textbf{Proof.} \\
\emph{1.} For $v_n<k<n-v_n$, according to the consistency of $\widehat{\Omega}_{n}(u_n)$ and the asymptotic normality of the MLE, we have as $n \rightarrow \infty$,
\begin{eqnarray*}
\left|\widehat{C}_{n,k}-C_{n,k}\right|&=&\frac{1}{q^{2}(\frac{k}{n})}.\frac{k^{2}(n-k)^{2}}{n^{3}}\left|\left(\widehat{\theta}(T_{1,k})-\widehat{\theta}(T_{k+1,n})\right)^{T}\left(\widehat{\Omega}_{n}(u_n)-\Omega\right)\left(\widehat{\theta}(T_{1,k})-\widehat{\theta}(T_{k+1,n})\right)\right|\\
&\leq&\frac{1}{q^{2}(\frac{k}{n})}.\frac{k^{2}(n-k)^{2}}{n^{3}}\left\|\widehat{\Omega}_{n}(u_n)-\Omega\right\| \left\|\widehat{\theta}(T_{1,k})-\widehat{\theta}(T_{k+1,n})\right\|^{2}\\
&\leq&C\frac{1}{q^{2}(\frac{k}{n})}.\frac{k(n-k)}{n^{2}}\left\|\widehat{\Omega}_{n}(u_n)-\Omega\right\| \Bigg[\left\|\sqrt{k}\left(\widehat{\theta}(T_{1,k})-\theta_0\right)\right\|^{2}\\
&~&~~~~~~~~~~~~~~~~~~~~~~~~~~~~~~~~~~~~~~~~~~~~~~~ +~\left\|\sqrt{n-k}\left(\widehat{\theta}(T_{k+1,n})-\theta_0\right)\right\|^{2}\Bigg]\\
&\leq&C\frac{1}{q^{2}(\frac{k}{n})}.\frac{k(n-k)}{n^{2}}o(1)O_P(1).
\end{eqnarray*}
Hence, as $n\rightarrow +\infty$,
\begin{eqnarray*}
\max_{v_n<k<n-v_n}\left|\widehat{C}_{n,k}-C_{n,k}\right|&\leq& o_P(1)\max_{v_n<k<n-v_n}\frac{1}{q^{2}(\frac{k}{n})}.\frac{k(n-k)}{n^{2}}\\
&\leq& o_P(1)\max_{\frac{v_n}{n}<\frac{k}{n}<1-\frac{v_n}{n}}\frac{1}{q^{2}(\frac{k}{n})}.\frac{k}{n}(1-\frac{k}{n})\\
&\leq& o_P(1)\sup_{0<\tau<1}\left(\frac{\sqrt{\tau(1-\tau)}}{q(\tau)}\right)^{2}.
%=o_P(1).
\end{eqnarray*}
The conclusion is easily obtained by using the fact of that $\sup_{0<\tau<1}\frac{\sqrt{\tau(1-\tau)}}{q(\tau)}<\infty$ if $I(q,c)<\infty$ for some $c>0$.\\
\\
\emph{2.} Under H$_0$, $(X_t,Y_t)_{t \in \Z}$ is a stationary and ergodic process, the same properties hold for
 $\left(\frac{\partial}{\partial \theta} \ell_t(\theta_0)\right)_{t \in \mathbb{Z}}$.
 Moreover, we have (see subsection \ref{Likelihood}),
\[\frac{\partial \ell_t(\theta)}{\partial\theta}=\left\{Y_{t}-A^{\prime}(\eta_t(\theta ))\right\}\frac{\partial \eta_t(\theta)}{\partial\theta}=\left\{Y_{t}-f^{t}_{\theta}\right\}\frac{\partial\left((A^{\prime})^{-1}(f^{t}_{\theta})\right)}{\partial\theta}.\]
Since $\frac{\partial \eta_t(\theta)}{\partial\theta}$ and $A^{\prime}(\eta_t(\theta ))$ are $\mathcal{F}_{t-1}$-measurable for any $\theta \in \Theta$, it holds that
$\mathbb{E}\left(\frac{\partial \ell_t(\theta_0)}{\partial\theta}|\mathcal{F}_{t-1}\right)=0$.\\
We also have,
\begin{eqnarray*}
\mathbb{E}\left(\left.\left(\frac{\partial \ell_t(\theta_0)}{\partial\theta}\right)\left(\frac{\partial \ell_t(\theta_0)}{\partial\theta}\right)^{T}\right|\mathcal{F}_{t-1}\right)&=&\mathbb{E}\left(\left.\left\{Y_{t}-A^{\prime}(\eta_t(\theta_0 ))\right\}^{2}\left(\frac{\partial \eta_t(\theta_0)}{\partial\theta}\right)\left(\frac{\partial \eta_t(\theta_0)}{\partial\theta}\right)^{T}\right|\mathcal{F}_{t-1}\right)\\
\\
&=&\mathbb{E}\left(\left.\left\{Y_{t}-A^{\prime}(\eta_t(\theta_0 ))\right\}^{2}\right|\mathcal{F}_{t-1}\right)\left(\frac{\partial \eta_t(\theta_0)}{\partial\theta}\right)\left(\frac{\partial \eta_t(\theta_0)}{\partial\theta}\right)^{T}\\
\\
&=&A^{\prime \prime}(\eta_t(\theta_0 ))\left(\frac{\partial \eta_t(\theta_0)}{\partial\theta}\right)\left(\frac{\partial \eta_t(\theta_0)}{\partial\theta}\right)^{T}.
\end{eqnarray*}
Therefore
\begin{eqnarray*}
\mathbb{E}\left[ \left(\frac{\partial \ell_t(\theta_0)}{\partial\theta}\right)\left(\frac{\partial \ell_t(\theta_0)}{\partial\theta}\right)^{T} \right]=\mathbb{E}\left(A^{\prime \prime}(\eta_t(\theta_0 ))\left(\frac{\partial \eta_t(\theta_0)}{\partial\theta}\right)\left(\frac{\partial \eta_t(\theta_0)}{\partial\theta}\right)^{T}\right)=\Omega.\\
\end{eqnarray*}
\emph{3.} We have,
\begin{eqnarray*}
\frac{\partial^{2}\ell_t(\theta_0)}{\partial \theta\partial \theta^{T}}&=& -A^{\prime\prime}(\eta_t(\theta_0 ))\left(\frac{\partial \eta_t(\theta_0)}{\partial\theta}\right)\left(\frac{\partial \eta_t(\theta_0)}{\partial\theta}\right)^{T}+\left\{Y_{t}-A^{\prime}(\eta_t(\theta_0 ))\right\}\frac{\partial^{2} \eta_t(\theta_0)}{\partial\theta\partial\theta^{T}}
\end{eqnarray*}
Moreover,
$$\mathbb{E}\left(\left.\left\{Y_{t}-A^{\prime}(\eta_t(\theta_0 ))\right\}\frac{\partial^{2} \eta_t(\theta_0)}{\partial\theta\partial\theta^{T}}\right|\mathcal{F}_{t-1}\right) = \mathbb{E}\left(\left.\left\{Y_{t}-A^{\prime}(\eta_t(\theta_0 ))\right\}\right|\mathcal{F}_{t-1}\right)\frac{\partial^{2} \eta_t(\theta_0)}{\partial\theta\partial\theta^{T}}=0.$$
Hence,
$$\mathbb{E}\left(\frac{\partial^{2}\ell_t(\theta_0)}{\partial \theta\partial \theta^{T}}\right)=-\mathbb{E}\left(A^{\prime\prime}(\eta_t(\theta_0 ))\left(\frac{\partial \eta_t(\theta_0)}{\partial\theta}\right)\left(\frac{\partial \eta_t(\theta_0)}{\partial\theta}\right)^{T}\right)=-\Omega.$$
~~\\
\emph{4.}  By applying $(\ref{G_n})$ with $\bar{\theta} = \widehat{\theta}_n(T_{1,n})$, it holds that
\begin{eqnarray*}
G_n(T_{1,n},\widehat{\theta}_n(T_{1,n}))= - \left(\frac{1}{n}\frac{\partial^{2}}{\partial \theta\partial \theta_i} L_n(T_{1,n},\theta_{n,i})\right)_{1\leq i \leq d}=-\frac{1}{n}\left(\sum_{t=1}^{n}\frac{\partial^{2}}{\partial \theta\partial \theta_i} \ell_t(\theta_{n,i})\right)_{1\leq i \leq d}
\end{eqnarray*}
where $\theta_{n,i}$ belongs between $\widehat{\theta}_n(T_{1,n})$ and $\theta_0$.
\\
For any $j=1,\cdots,d$, we have
\begin{align}
\nonumber \left|\frac{1}{n}\sum_{t=1}^{n}\frac{\partial^{2}}{\partial \theta_j\partial \theta_i} \ell_t(\theta_{n,i})-\mathbb{E}\left(\frac{\partial^{2}}{\partial \theta_j\partial \theta_i} \ell_0(\theta_{0})\right)\right|
 \nonumber  &\leq \left|\frac{1}{n}\sum_{t=1}^{n}\frac{\partial^{2}}{\partial \theta_j\partial \theta_i} \ell_t(\theta_{n,i})-\mathbb{E}\left(\frac{\partial^{2}}{\partial \theta_j\partial \theta_i} \ell_0(\theta_{n,i})\right)\right|\\
 \nonumber &+\left|\mathbb{E}\left(\frac{\partial^{2}}{\partial \theta_j\partial \theta_i} \ell_0(\theta_{n,i})\right)-\mathbb{E}\left(\frac{\partial^{2}}{\partial \theta_j\partial \theta_i} \ell_0(\theta_{0})\right)\right| \\
 \nonumber &\leq \left\|\frac{1}{n}\sum_{t=1}^{n}\frac{\partial^{2}}{\partial \theta_j\partial \theta_i} \ell_t(\theta)-\mathbb{E}\left(\frac{\partial^{2}}{\partial \theta_j\partial \theta_i} \ell_0(\theta)\right)\right\|_\Theta\\
& + \left|\mathbb{E}\left(\frac{\partial^{2}}{\partial \theta_j\partial \theta_i} \ell_0(\theta_{n,i})\right)-\mathbb{E}\left(\frac{\partial^{2}}{\partial \theta_j\partial \theta_i} \ell_0(\theta_{0})\right)\right|  \label{lem1_4_G}
% &~&~\stackrel{a.s}{\longrightarrow}0~~\textrm{as}~n\rightarrow +\infty.
\end{align}
 We have
 \[
 \left|\mathbb{E}\left(\frac{\partial^{2}}{\partial \theta_j\partial \theta_i} \ell_0(\theta_{n,i})\right)-\mathbb{E}\left(\frac{\partial^{2}}{\partial \theta_j\partial \theta_i} \ell_0(\theta_{0})\right)\right| \limiten  0.
\]
 Moreover, the sequence $ \big(\frac{\partial^2}{\partial \theta_i \partial \theta_j}\ell_t(\theta) \big)_{t \in \Z}$ is stationary ergodicity and we have
 \[
 \frac{\partial^2}{\partial \theta_i \partial \theta_j}\ell_t(\theta) = -A^{\prime \prime}(\eta_t(\theta ))\left(\frac{\partial}{\partial\theta_i}\eta_t(\theta)\right) \left(\frac{\partial}{\partial\theta_j}\eta_t(\theta)\right) +  \big\{Y_{t}-A^{\prime}(\eta_t(\theta ))\big\}\frac{\partial^2}{\partial\theta_i\partial\theta_j}\eta_t(\theta)
\]
which shows from (\textbf{A6}) that
 $ \E [ \| \frac{\partial^2}{\partial \theta_i \partial \theta_j}\ell_t(\theta) \|_{\Theta} ]  < + \infty $.
 Hence, by the uniform strong law of large numbers, it holds that
\[\left\|\frac{1}{n}\sum_{t=1}^{n}\frac{\partial^{2}}{\partial \theta_i\partial \theta_j} \ell_t(\theta)-\mathbb{E}\left(\frac{\partial^{2}}{\partial \theta_i\partial \theta_j} \ell_0(\theta)\right)\right\|_\Theta \underset{n \to\infty}{\longrightarrow} 0.\]
Thus, according to (\ref{lem1_4_G}), we get
\[  \left|\frac{1}{n}\sum_{t=1}^{n}\frac{\partial^{2}}{\partial \theta_j\partial \theta_i} \ell_t(\theta_{n,i})-\mathbb{E}\left(\frac{\partial^{2}}{\partial \theta_j\partial \theta_i} \ell_0(\theta_{0})\right)\right|
  \limiten  0 \]
which shows that
\[
G_n(T_{1,n},\widehat{\theta}_n(T_{1,n}))=-\frac{1}{n}\left(\sum_{t=1}^{n}\frac{\partial^{2}}{\partial \theta\partial \theta_i} \ell_t(\theta_{n,i})\right)_{1\leq i \leq d} \limitepsn   -\mathbb{E}\left(\frac{\partial^{2}}{\partial \theta\partial \theta^{T}} \ell_0(\theta_{0})\right)=\Omega.
\]
  By going along similar lines, one can prove that
   $ \left.\dfrac{1}{n}\sum_{t=1}^{n}\left({A}^{\prime \prime}\left( \eta_t(\theta)\right)\left(\frac{\partial \eta_t(\theta)}{\partial \theta}\right)\left(\frac{\partial \eta_t(\theta)}{\partial \theta}\right)^{T}\right)\right|_{\theta=\widehat{\theta}_n }
     \stackrel{a.s}{\longrightarrow}\Omega$ as $n\rightarrow +\infty$. \\
 This completes the proof of the lemma.

 ~ \\

Now, let us use the Lemma \ref{Lem_H0} to show that
\[C_n\stackrel{\mathcal{D}}{\longrightarrow}\sup_{0\leq\tau\leq1}\frac{\left\|W_d(\tau)\right\|^{2}}{q^{2}(\tau)}~~ \textrm{as}~~n\rightarrow \infty.\]
Let $v_n \leq k \leq n-v_n$. By applying (\ref{Eq_Taylor}) with $\bar{\theta}=\widehat{\theta}_n(T_{1,k})$ and $T_{k,k^\prime}=T_{1,k}$, we get
\begin{eqnarray}\label{Eq_T_1,k}
G_n(T_{1,k},\widehat{\theta}_n(T_{1,k})) \cdot (\widehat{\theta}_n(T_{1,k})-\theta_0)=\frac{1}{k}\left(\frac{\partial}{\partial \theta} L_n(T_{1,k},\theta_0)-\frac{\partial}{\partial \theta} L_n(T_{1,k},\widehat{\theta}_n(T_{1,k}))\right).
\end{eqnarray}
With $\bar{\theta}=\widehat{\theta}_n(T_{k+1,n})$ and $T_{k,k^\prime}=T_{k+1,n}$, (\ref{Eq_Taylor}) becomes
\begin{eqnarray}\label{Eq_T_k+1,n}
G_n(T_{k+1,n},\widehat{\theta}_n(T_{k+1,n})) \cdot (\widehat{\theta}_n(T_{k+1,n})-\theta_0)=\frac{1}{n-k}\left(\frac{\partial}{\partial \theta} L_n(T_{k+1,n},\theta_0)-\frac{\partial}{\partial \theta} L_n(T_{k+1,n},\widehat{\theta}_n(T_{k+1,n}))\right).
\end{eqnarray}
As $n\rightarrow +\infty$, we have
\begin{eqnarray*}
&~&\left\|G_n(T_{1,k},\widehat{\theta}_n(T_{1,k}))-\Omega\right\|=\left\|G_n(T_{k+1,n},\widehat{\theta}_n(T_{k+1,n}))-\Omega\right\|=o(1),\\
\\ &~&\sqrt{k}\left(\widehat{\theta}_n(T_{1,k})-\theta_0\right)=O_P(1)~\textrm{and}~\sqrt{n-k}\left(\widehat{\theta}_n(T_{k+1,n})-\theta_0\right)=O_P(1).
\end{eqnarray*}
According (\ref{Eq_T_1,k}), we have
\begin{eqnarray*}
\sqrt{k}\Omega\left(\widehat{\theta}_n(T_{1,k})-\theta_0\right)&=&\frac{1}{\sqrt{k}}\left(\frac{\partial}{\partial \theta} L_n(T_{1,k},\theta_0)-\frac{\partial}{\partial \theta} L_n(T_{1,k},\widehat{\theta}_n(T_{1,k}))\right)\\
\\
&~&~-\sqrt{k}\left(\left(G_n(T_{1,k},\widehat{\theta}_n(T_{1,k}))-\Omega\right)\left(\widehat{\theta}_n(T_{1,k})-\theta_0\right)\right)\\
&=&\frac{1}{\sqrt{k}}\left(\frac{\partial}{\partial \theta} L_n(T_{1,k},\theta_0)-\frac{\partial}{\partial \theta} L_n(T_{1,k},\widehat{\theta}_n(T_{1,k}))\right) +o_P(1).
\end{eqnarray*}
It is equivalent to
\begin{eqnarray}\label{Eq_a}
\Omega\left(\widehat{\theta}_n(T_{1,k})-\theta_0\right)=\frac{1}{k}\left(\frac{\partial}{\partial \theta} L_n(T_{1,k},\theta_0)-\frac{\partial}{\partial \theta} L_n(T_{1,k},\widehat{\theta}_n(T_{1,k}))\right) +o_P\left(\frac{1}{\sqrt{k}}\right).
\end{eqnarray}
 For $n$ large enough, $\widehat{\theta}_n(T_{1,k})$ is an interior point of $\Theta$ and we have $\frac{\partial}{\partial \theta}L_n(T_{1,k},\widehat{\theta}_n(T_{1,k})) = 0$.
Hence, for $n$ large enough, we get from (\ref{Eq_a})
\begin{eqnarray*}
\Omega\left(\widehat{\theta}_n(T_{1,k})-\theta_0\right)
 &=&\frac{1}{k}\frac{\partial}{\partial \theta} L_n(T_{1,k},\theta_0)+o_P\left(\frac{1}{\sqrt{k}}\right).
\end{eqnarray*}
Similarly, we can use (\ref{Eq_T_k+1,n}) to obtain
$$\Omega\left(\widehat{\theta}_n(T_{k+1,n})-\theta_0\right)=\frac{1}{n-k}\frac{\partial}{\partial \theta} L_n(T_{k+1,n},\theta_0)+o_P\left(\frac{1}{\sqrt{n-k}}\right).$$
The subtraction of the two above equations gives
\begin{eqnarray*}
\Omega\left(\widehat{\theta}_n(T_{1,k})-\widehat{\theta}_n(T_{k+1,n})\right)&=&\frac{1}{k}\frac{\partial}{\partial \theta} L_n(T_{1,k},\theta_0)-\frac{1}{n-k}\frac{\partial}{\partial \theta} L_n(T_{k+1,n},\theta_0)+o_P\left(\frac{1}{\sqrt{k}}+\frac{1}{\sqrt{n-k}}\right)\\
&=&\frac{1}{k}\frac{\partial}{\partial \theta} L_n(T_{1,k},\theta_0)-\frac{1}{n-k}\left(\frac{\partial}{\partial \theta} L_n(T_{1,n},\theta_0)-\frac{\partial}{\partial \theta} L_n(T_{1,k},\theta_0)\right)\\
&~&~+o_P\left(\frac{1}{\sqrt{k}}+\frac{1}{\sqrt{n-k}}\right)\\
\\
&=&\frac{n}{k(n-k)}\left(\frac{\partial}{\partial \theta} L_n(T_{1,k},\theta_0)-\frac{n}{k}.\frac{\partial}{\partial \theta} L_n(T_{1,n},\theta_0)\right)+o_P\left(\frac{1}{\sqrt{k}}+\frac{1}{\sqrt{n-k}}\right).
\end{eqnarray*}
i.e.
\begin{eqnarray*}
\frac{k(n-k)}{n^{\frac{3}{2}}}\Omega(\widehat{\theta}_n(T_{1,k})-\widehat{\theta}_n(T_{k+1,n}))&=&\frac{1}{\sqrt{n}}\left(\frac{\partial}{\partial \theta} L_n(T_{1,k},\theta_0)-\frac{k}{n}.\frac{\partial}{\partial \theta} L_n(T_{1,n},\theta_0)\right)\\
&~&~+o_P\left(\frac{\sqrt{k(n-k)}}{n} + \frac{ \sqrt{n-k}}{\sqrt{n}}\right)\\
\\
&=&\frac{1}{\sqrt{n}}\left(\frac{\partial}{\partial \theta} L_n(T_{1,k},\theta_0)-\frac{k}{n}.\frac{\partial}{\partial \theta} L_n(T_{1,n},\theta_0)\right) +o_P(1).
\end{eqnarray*}
According the above equation, we have
\begin{eqnarray}\label{Eq_b}
\frac{k(n-k)}{n^{\frac{3}{2}}}\Omega^{1/2}(\widehat{\theta}_n(T_{1,k})-\widehat{\theta}_n(T_{k+1,n}))=\frac{\Omega^{-1/2}}{\sqrt{n}}\left(\frac{\partial}{\partial \theta} L_n(T_{1,k},\theta_0)-\frac{k}{n}.\frac{\partial}{\partial \theta} L_n(T_{1,n},\theta_0)\right) +o_P(1).
\end{eqnarray}
Recall that for any $0<\tau<1$  $$\frac{\partial}{\partial \theta}L_n(T_{1,[n \tau]},\theta_0)=\sum_{t=1}^{[n \tau]}\frac{\partial}{\partial \theta}\ell_t(\theta_0)$$
where $[n \tau]$ is the integer part of $n \tau$.\\
The process $\left(\frac{\partial}{\partial \theta} \ell_t(\theta_0),\mathcal{F}_{t}\right)_{t \in \mathbb{Z}}$ is a stationary ergodic square integrable martingale difference process with covariance matrix $\Omega$ (see Lemma \ref{Lem_H0}). By applying the Central limit theorem for the martingale difference sequence (see Billingsley (1968)), we have
\begin{eqnarray*}
\frac{1}{\sqrt{n}}\left(\frac{\partial}{\partial \theta}L_n(T_{1,[n \tau]},\theta_0)-\frac{[n \tau]}{n}\frac{\partial}{\partial \theta}L_n(T_{1,n },\theta_0)\right)&=&\frac{1}{\sqrt{n}}\left(\sum_{t=1}^{[n \tau]}\frac{\partial}{\partial \theta}\ell_t(\theta_0)-\frac{[n \tau]}{n}\sum_{t=1}^{n }\frac{\partial}{\partial \theta}\ell_t(\theta_0)\right)\\
\\
&~&~\stackrel{\mathcal{D}}{\longrightarrow} B_{\Omega}(\tau)-\tau B_{\Omega}(1)~~~\textrm{as}~~n\rightarrow +\infty.
\end{eqnarray*}
where $B_{\Omega}$ is a Gaussian process with covariance matrix $\min(s,t)\Omega$.\\
Hence,
$$\frac{1}{\sqrt{n}}\Omega^{-1/2}\left(\frac{\partial}{\partial \theta}L_n(T_{1,[n \tau]},\theta_0)-\frac{[n \tau]}{n}\frac{\partial}{\partial \theta}L_n(T_{1,n },\theta_0)\right)\stackrel{\mathcal{D}}{\longrightarrow} B_{d}(\tau)-\tau B_{d}(1)=W_d(\tau)~~~\textrm{as}~~n\rightarrow +\infty$$
in $D\left(\left[0,1\right]\right)$ where $B_d$ is a $d$-dimensional standard motion, and $W_d$ is a $d$-dimensional Brownian bridge.\\
From (\ref{Eq_b}), as $n\rightarrow +\infty$, we have
\begin{eqnarray*}
C_{n,\left[n\tau\right]}&=&\frac{1}{q^{2}\left(\frac{[n\tau]}{n}\right)}\frac{[n\tau]^{2}(n-[n\tau])^{2}}{n^{3}}\left(\widehat{\theta}(T_{1,[n\tau]})-\widehat{\theta}(T_{[n\tau]+1,n})\right)^{T}\Omega \left(\widehat{\theta}(T_{1,[n\tau]})-\widehat{\theta}(T_{[n\tau]+1,n})\right)\\
&=&\frac{1}{q^{2}\left(\frac{[n\tau]}{n}\right)}\left\|\frac{[n\tau](n-[n\tau])}{n^{3/2}}\Omega^{1/2}\left(\widehat{\theta}(T_{1,[n\tau]})-\widehat{\theta}(T_{[n\tau]+1,n})\right)\right\|^{2}\\
&=&\frac{1}{q^{2}\left(\frac{[n\tau]}{n}\right)}\left\|\frac{\Omega^{-1/2}}{\sqrt{n}}\left(\frac{\partial}{\partial \theta} L_n(T_{1,[n\tau]},\theta_0)-\frac{[n\tau]}{n}.\frac{\partial}{\partial \theta} L_n(T_{1,n},\theta_0)\right) +o_P(1)\right\|^{2}\\
&=&\frac{1}{q^{2}\left(\frac{[n\tau]}{n}\right)}\left\|\frac{\Omega^{-1/2}}{\sqrt{n}}\left(\sum_{t=1}^{[n \tau]}\frac{\partial}{\partial \theta}\ell_t(\theta_0)-\frac{[n \tau]}{n}\sum_{t=1}^{n }\frac{\partial}{\partial \theta}\ell_t(\theta_0)\right) \right\|^{2} +o_P(1)\\
\\
&~&~~\stackrel{\mathcal{D}}{\longrightarrow} \frac{\left\|W_d(\tau)\right\|^{2}}{q^{2}(\tau)}~\textrm{in}~ \mathcal{D}\left(\left[0,1\right]\right).
\end{eqnarray*}
According to the properties of $q$, we have for any $0<\epsilon<1/2$
\begin{eqnarray*}
\max_{[n\epsilon]<k<n-[n\epsilon]}C_{n,k}&=&\max_{[n\epsilon]<k<n-[n\epsilon]}\frac{1}{q^{2}\left(\frac{k}{n}\right)}\frac{k^{2}(n-k)^{2}}{n^{3}}\left(\widehat{\theta}(T_{1,k})-\widehat{\theta}(T_{k+1,n})\right)^{T}\Omega \left(\widehat{\theta}(T_{1,k})-\widehat{\theta}(T_{k+1,n})\right)\\
\\
&=&\sup_{\epsilon<\tau<1-\epsilon}\frac{1}{q^{2}\left(\frac{[n\tau]}{n}\right)}\frac{[n\tau]^{2}(n-[n\tau])^{2}}{n^{3}}\left(\widehat{\theta}(T_{1,[n\tau]})-\widehat{\theta}(T_{[n\tau]+1,n})\right)^{T}\Omega\\ &~&~~~\times\left(\widehat{\theta}(T_{1,[n\tau]})-\widehat{\theta}(T_{[n\tau]+1,n})\right)\\
\\
&=&\sup_{\epsilon<\tau<1-\epsilon}\frac{1}{q^{2}\left(\frac{[n\tau]}{n}\right)}\left\|\frac{[n\tau](n-[n\tau])}{n^{3/2}}\Omega^{1/2}\left(\widehat{\theta}(T_{1,[n\tau]})-\widehat{\theta}(T_{[n\tau]+1,n})\right)\right\|^{2}\\
&=&\sup_{\epsilon<\tau<1-\epsilon}\frac{1}{q^{2}\left(\frac{[n\tau]}{n}\right)}\left\|\frac{\Omega^{-1/2}}{\sqrt{n}}\left(\frac{\partial}{\partial \theta} L_n(T_{1,[n\tau]},\theta_0)-\frac{[n\tau]}{n}.\frac{\partial}{\partial \theta} L_n(T_{1,n},\theta_0)\right) +o_P(1)\right\|^{2}\\
&=&\sup_{\epsilon<\tau<1-\epsilon}\frac{1}{q^{2}\left(\frac{[n\tau]}{n}\right)}\left\|\frac{\Omega^{-1/2}}{\sqrt{n}}\left(\sum_{t=1}^{[n \tau]}\frac{\partial}{\partial \theta}\ell_t(\theta_0)-\frac{[n \tau]}{\sqrt{n}}\sum_{t=1}^{n }\frac{\partial}{\partial \theta}\ell_t(\theta_0)\right) \right\|^{2} +o_P(1)\\
\\
&~&~~\stackrel{\mathcal{D}}{\longrightarrow} \sup_{\epsilon<\tau<1-\epsilon}\frac{\left\|W_d(\tau)\right\|^{2}}{q^{2}(\tau)}~\textrm{in}~ \mathcal{D}\left(\left[0,1\right]\right).
\end{eqnarray*}
Hence, as $n\rightarrow +\infty$, we have shown that
$$C_{n,\left[n\tau\right]}\stackrel{\mathcal{D}}{\longrightarrow} \frac{\left\|W_d(\tau)\right\|^{2}}{q^{2}(\tau)}~\textrm{in}~ \mathcal{D}\left(\left[0,1\right]\right)$$
and for all $0<\epsilon<1/2$
$$\max_{[n\epsilon]<k<n-[n\epsilon]}C_{n,k}=\sup_{\epsilon<\tau<1-\epsilon}C_{n,[n\tau]}\stackrel{\mathcal{D}}{\longrightarrow} \sup_{\epsilon<\tau<1-\epsilon}\frac{\left\|W_d(\tau)\right\|^{2}}{q^{2}(\tau)}.$$
In addition, since $I(q,c)<+\infty$ for some $c>0$, one can show that almost surely $$\lim_{\tau\rightarrow 0}\frac{\left\|W_d(\tau)\right\|}{q(\tau)}<\infty~~\textrm{and}~~\lim_{\tau\rightarrow 1}\frac{\left\|W_d(\tau)\right\|}{q(\tau)}<\infty.$$
Hence, for $n$ large enough we have
$$C_n=\max_{v_n<k<n-v_n}C_{n,k}=\sup_{\frac{v_n}{n}<\tau<1-\frac{v_n}{n}}C_{n,[n\tau]}\stackrel{\mathcal{D}}{\longrightarrow}\sup_{0<\tau<1}\frac{\left\|W_d(\tau)\right\|^{2}}{q^{2}(\tau)}~~\textrm{as}~~ n\rightarrow +\infty.~~~~~~~~~~\hfill\square\\$$

\subsection{Proof of Theorem \ref{Thm_Asy_H1}}
We assume that the trajectory $(Y_1,\cdots,Y_n)$ satisfies
\begin{equation} \label{Eq_H1}
Y_{t}=\left\{
\begin{array}{ll}
Y^{(1)}_t~~\textrm{for}~~t\leq t^{*},\\
\\
Y^{(2)}_t~~\textrm{for}~~t>t^{*},\\
\end{array}
\right.
\end{equation} \\
where $t^{*}=[\tau n]$ with $0<\tau<1$ and $\{Y^{(i)}_t, t \in \mathbb{Z}\}$ ($i=1,2$) is a stationary solution of the model (\ref{Model}) depending on $\theta^{*}_i$ with $\theta^{*}_1 \neq \theta^{*}_2$.\\
\\
We have
$$\widehat{C}_{n,t^{*}}=\frac{1}{q^{2}(\frac{t^{*}}{n})}\frac{{t^{*}}^{2}(n-t^{*})^{2}}{n^{3}}\left(\widehat{\theta}(T_{1,t^{*}})-\widehat{\theta}(T_{t^{*}+1,n})\right)^{T}\widehat{\Omega}_{n}(u_n)\left(\widehat{\theta}(T_{1,t^{*}})-\widehat{\theta}(T_{t^{*}+1,n})\right)$$
and $\widehat{C}_{n}=\max_{v_n\leq k \leq n-v_n}\widehat{C}_{n,k}\geq \widehat{C}_{n,t^{*}}$. Then, to prove the Theorem \ref{Thm_Asy_H1}, it suffices to show that $\widehat{C}_{n,t^{*}}\stackrel{P}{\longrightarrow} +\infty$ as $n\rightarrow +\infty $.\\
The likelihood function of the stationary process $\{Y^{(1)}_t, t \in \mathbb{Z}\}$ computed on $\left\{1,2,\cdots,t^{*}\right\}$ is given by $$L_n({T_{1,t^{*}}},\theta)=\sum_{t=1}^{t^{*}}\ell_{t,1}(\theta)$$
where  $\ell_{t,1}(\theta)=\{\widehat{\eta}_{t,1}(\theta )Y^{(1)}_{t}-A(\widehat{\eta}_{t,1}(\theta ))\}$ with $\widehat{\eta}_{t,1}=(A^{\prime})^{-1}(f^{t,1}_\theta)$ and $f^{t,1}_\theta:=f^{t}_\theta(Y^{(1)}_{t-1},Y^{(1)}_{t-2},\cdots)$. \\
Recall that the matrix used to construct the test statistic is $\widehat{\Omega}_{n}(u_n)$ given by
\begin{eqnarray*}
\widehat{\Omega}_{n}(u_n)&=&\frac{1}{2}\Bigg[\left.\frac{1}{u_n}\sum_{t=1}^{u_n}{A}^{\prime \prime}\left(\widehat{\eta}_t(\theta)\right)\left(\frac{\partial \widehat{\eta}_t(\theta)}{\partial \theta}\right)\left(\frac{\partial \widehat{\eta}_t(\theta)}{\partial \theta}\right)^{T}\right|_{\theta=\widehat{\theta}_n(T_{1,u_n})}\\
\\
\\
 &~&~~+~\left.\frac{1}{n-u_n}\sum_{t=u_n+1}^{n}{A}^{\prime \prime}\left(\widehat{\eta}_t(\theta)\right)\left(\frac{\partial \widehat{\eta}_t(\theta)}{\partial \theta}\right)\left(\frac{\partial \widehat{\eta}_t(\theta)}{\partial \theta}\right)^{T}\right|_{\theta= \widehat{\theta}_n(T_{u_n+1,n})}\Bigg].\\
 \end{eqnarray*}
According to the asymptotic proprieties of the maximum likelihood estimator, we have $$\widehat{\theta}_n(T_{1,t^{*}})\stackrel{a.s}{\longrightarrow} \theta^{*}_1,~~~\widehat{\theta}_n(T_{1,u_n})\stackrel{a.s}{\longrightarrow} \theta^{*}_1~~\textrm{as}~ n\longrightarrow +\infty~~~~~~~~~~~~~~~~~~~~~~~~~~~~~~~~$$
and
$$\frac{1}{u_n}\sum_{t=1}^{u_n} {A}^{\prime \prime}\left(\widehat{\eta}_t(\theta)\right)\left(\frac{\partial \widehat{\eta}_t(\theta)}{\partial \theta}\right)\left(\frac{\partial \widehat{\eta}_t(\theta)}{\partial \theta}\right)^{T} \stackrel{a.s}{\longrightarrow}\Omega^{(1)} ~\textrm{as}~ n\longrightarrow +\infty$$
where
$$\Omega^{(1)}=\mathbb{E}\left[{A}^{\prime \prime}\left(\widehat{\eta}_{0,1}(\theta^{*}_1)\right)\left(\frac{\partial \widehat{\eta}_{0,1}(\theta^{*}_1)}{\partial \theta}\right)\left(\frac{\partial \widehat{\eta}_{0,1}(\theta^{*}_1)}{\partial \theta}\right)^{T}\right].$$
~~\\
Define $L_{n,2}({T_{t^{*}+1,n}},\theta)=\sum_{t=t^{*}+1}^{n}\ell_{t,2}(\theta)$ the likelihood of the stationary process $\{Y^{(2)}_t, t \in \mathbb{Z}\}$ computed on $\left\{t^{*}+1,\cdots,n\right\}$, where $$\ell_{t,2}(\theta)=\{\widehat{\eta}_{t}(\theta )Y_{t}-A(\widehat{\eta}_{t}(\theta ))\}=\{\widehat{\eta}_{t,2}(\theta )Y^{(2)}_{t}-A(\widehat{\eta}_{t,2}(\theta ))\}$$
with $$\widehat{\eta}_{t,2}=(A^{\prime})^{-1}(f^{t,2}_\theta) ~~\textrm{and}~~f^{t,2}_\theta:=f^{t}_\theta(Y^{(2)}_{t-1},Y^{(2)}_{t-2},\cdots).$$
The asymptotic proprieties imply that
\begin{eqnarray*}
\widehat{\theta}_n(T_{t^{*}+1,n})=\textrm{argmax}_{\theta \in \Theta}L_{n}({T_{t^{*}+1,n}},\theta) = \textrm{argmax}_{\theta \in \Theta}L_{n,2}({T_{t^{*}+1,n}},\theta)&\stackrel{a.s}{\longrightarrow}& \theta^{*}_2 ~~~\textrm{as}~~~ n\longrightarrow +\infty.\\
\end{eqnarray*}
Recall that, by definition, the two matrices in the formula of $\widehat{\Omega}_n(u_n)$ are positive semi-definite  and the first one converges a.s. to $\Omega^{(1)}$ which is positive definite.\\
For $n$ large enough, we can write
\begin{eqnarray*}
\widehat{C}_{n}&\geq& \widehat{C}_{n,t^{*}}\\
&\geq& \frac{1}{q^{2}(\frac{t^{*}}{n})}\frac{{t^{*}}^{2}(n-t^{*})^{2}}{n^{3}}\left(\widehat{\theta}(T_{1,t^{*}})-\widehat{\theta}(T_{t^{*}+1,n})\right)^{T}\\
\\
 &~&~~\times\frac{1}{2}\Bigg[\left.\frac{1}{u_n}\sum_{t=1}^{u_n}{A}^{\prime \prime}\left(\widehat{\eta}_t(\theta)\right)\left(\frac{\partial \widehat{\eta}_t(\theta)}{\partial \theta}\right)\left(\frac{\partial \widehat{\eta}_t(\theta)}{\partial \theta}\right)^{T}\right|_{\theta=\widehat{\theta}_n(T_{1,u_n})}\\
\\
 &~&~~~~~~~~~+~\left.\frac{1}{n-u_n}\sum_{t=u_n+1}^{n}{A}^{\prime \prime}\left(\widehat{\eta}_t(\theta)\right)\left(\frac{\partial \widehat{\eta}_t(\theta)}{\partial \theta}\right)\left(\frac{\partial \widehat{\eta}_t(\theta)}{\partial \theta}\right)^{T}\right|_{\theta= \widehat{\theta}_n(T_{u_n+1,n})}\!\Bigg]\\
&~&~~\times\left(\widehat{\theta}(T_{1,t^{*}})-\widehat{\theta}(T_{t^{*}+1,n})\right)\\
\\
&\geq&\frac{1}{q^{2}(\frac{t^{*}}{n})}\frac{{t^{*}}^{2}(n-t^{*})^{2}}{n^{3}}\left(\widehat{\theta}(T_{1,t^{*}})-\widehat{\theta}(T_{t^{*}+1,n})\right)^{T}\times~\\
\\
 &~&~~\frac{1}{2}\left[\left.\frac{1}{u_n}\sum_{t=1}^{u_n}{A}^{\prime \prime}\left(\widehat{\eta}_t(\theta)\right)\left(\frac{\partial \widehat{\eta}_t(\theta)}{\partial \theta}\right)\left(\frac{\partial \widehat{\eta}_t(\theta)}{\partial \theta}\right)^{T}\right|_{\theta=\widehat{\theta}_n(T_{1,u_n})}\right]\times\left(\widehat{\theta}(T_{1,t^{*}})-\widehat{\theta}(T_{t^{*}+1,n})\right)\\
 \\
\\ &\geq&\frac{1}{\sup_{0<\tau<t}q^{2}(\tau)}n\left(\tau^{*}(1-\tau^{*})\right)^{2}\left(\widehat{\theta}(T_{1,t^{*}})-\widehat{\theta}(T_{t^{*}+1,n})\right)^{T}\times~\\
\\
 &~&~~\frac{1}{2}\left[\left.\frac{1}{u_n}\sum_{t=1}^{u_n}{A}^{\prime \prime}\left(\widehat{\eta}_t(\theta)\right)\left(\frac{\partial \widehat{\eta}_t(\theta)}{\partial \theta}\right)\left(\frac{\partial \widehat{\eta}_t(\theta)}{\partial \theta}\right)^{T}\right|_{\theta=\widehat{\theta}_n(T_{1,u_n})}\right]\times\left(\widehat{\theta}(T_{1,t^{*}})-\widehat{\theta}(T_{t^{*}+1,n})\right)\\
\\
\\
&\geq& C \times n \left(\widehat{\theta}(T_{1,t^{*}})-\widehat{\theta}(T_{t^{*}+1,n})\right)^{T}\times\frac{1}{2}\left[\left.\frac{1}{u_n}\sum_{t=1}^{u_n}{A}^{\prime \prime}\left(\widehat{\eta}_t(\theta)\right)\left(\frac{\partial \widehat{\eta}_t(\theta)}{\partial \theta}\right)\left(\frac{\partial \widehat{\eta}_t(\theta)}{\partial \theta}\right)^{T}\right|_{\theta=\widehat{\theta}_n(T_{1,u_n})}\right]\\
\\
&~&~\times~\left(\widehat{\theta}(T_{1,t^{*}})-\widehat{\theta}(T_{t^{*}+1,n})\right).
\end{eqnarray*}
Moreover, as $n \rightarrow+\infty$, we have $$
\widehat{\theta}(T_{1,t^{*}})-\widehat{\theta}(T_{t^{*}+1,n})\stackrel{a.s}{\longrightarrow}\theta^{*}_{1}-\theta^{*}_{2}\neq 0
~~~~\textrm{and} ~~~~\frac{1}{u_n}\sum_{t=1}^{u_n}{A}^{\prime \prime}\left(\widehat{\eta}_t(\theta)\right)\left(\frac{\partial \widehat{\eta}_t(\theta)}{\partial \theta}\right)\left(\frac{\partial \widehat{\eta}_t(\theta)}{\partial \theta}\right)^{T}\stackrel{a.s}{\longrightarrow}\Omega^{(1)}.
$$
Hence, $\widehat{C}_{n}\stackrel{a.s}{\longrightarrow}+\infty~~\textrm{as}~~ n \rightarrow+\infty.~~~~~~\hfill\square\\$

\end{document}